\documentclass[a4paper,12pt]{article} 
\newdimen\paperhight
\usepackage{amsmath}
\usepackage{amssymb}
\usepackage{amsfonts}
\usepackage[usenames]{color}

\newcommand{\ol}{\overline}

\newcommand{\pd}{\partial} 

\newcommand{\pr}{\par \vspace{3mm}\noindent [{\bf Proof}] \qquad}
\newcommand{\prend}{\hfill \qed \par \vspace{3mm}}

\newcommand{\qed}{\quad\hbox{\rule[-2pt]{3pt}{6pt}}\par\vspace{3mm}}

\newcommand{\1}{{\bf 1}} 
 
\newcommand{\C}{\mathbb C} 
\newcommand{\Z}{\mathbb Z} 
\newcommand{\Q}{\mathbb Q} 
 
\newcommand{\N}{\mathbb N} 
\newcommand{\R}{\mathbb R}

\newcommand{\CM}{{\cal M}}

\newcommand{\CV}{{\cal V}}

\newcommand{\CY}{{\cal Y}}

\newcommand{\be}{\beta}

\newcommand{\la}{\lambda}

\newcommand{\wt}{{\rm wt}}
\newcommand{\Hom}{{\rm Hom}}

\newcommand{\soc}{{\rm soc}}
\newcommand{\End}{{\rm End}}
\newcommand{\Ker}{{\rm Ker}}

\newtheorem{thm}{Theorem}[section]
\newtheorem{prn}[thm]{Proposition}
\newtheorem{dfn}[thm]{Definition}
\newtheorem{cry}[thm]{Corollary}
\newtheorem{lmm}[thm]{Lemma}

\newtheorem{rmk}{Remark}

\begin{document}
\title{A theory of tensor products for vertex operator algebra 
satisfying $C_2$-cofiniteness} 
\author{Masahiko Miyamoto\footnote{Supoorted by the Grants-in-Aids 
for Scientific Research, No. 1344002, The Ministry of Education, 
Science and Culture, Japan}}
\date{\begin{tabular}{c}
Institute of Mathematics, \cr
University of Tsukuba \cr
Tsukuba, 305 Japan\end{tabular}}
\maketitle

\begin{abstract}
The recent researchs show that $C_2$-cofiniteness is a natural 
conditition to consider a vertex operator algebra with finitely 
many simple modules. Therefore, 
we extended the tensor product theory of vertex operator algebras 
developed by Huang and Lepowsky without assuming the compatibility 
condition nor the semisimplicity of grading operator so that we could 
apply it to all vertex operator algebras satisfying only 
$C_2$-cofiniteness. 
We also showed that the tensor product theory develops naturally 
if we include not only ordinary modules, but also weak modules with 
a composition series of finite length (we call it an Artin module). 
Actually, a $C_2$-cofiniteness on $V$ is enough to show that a 
tensor product of two Artin modules is again an Arting module and we 
have natural commutativity and associativity of 
tensor products. Namely, the category of Artin modules becomes a 
braided tensor category.

Our main purpose is an application of the tensor product theory under 
$C_2$-cofiniteness. We determined the representation theory of orbifold 
models. For example, if a vertex operator algebra $V$ has a finite 
automorphism group $G$ and the fixed point vertex operator subalgebra 
$V^G$ is $C_2$-cofinite, then for $g\in G$ and any irreducible 
$V^{\langle g\rangle }$-module $W$, there is an element 
$h\in \langle g\rangle $ such that 
$W$ is contained in some $h$-twisted $V$-module. Furthermore, if 
$V^G$ is rational, then $V^{\langle g\rangle }$ is also rational 
for any $g\in G$. 
\end{abstract}

\section{Introduction}
The tensor product theory for the category of modules for a 
vertex operator algebra was developed by a sequential studies of 
Huang and Lepowsky under the assumptions of the compatibility 
condition and semisimplicity of grading operator, 
\cite{HL1}-\cite{HL3}\cite{H1},\cite{H2}. Above all, 
Huang showed that if every irreducible module 
$W$ satisfies $C_1$-cofiniteness ($\dim W/C_1(W)<\infty$), then 
$n$-point correlation functions defined by products of intertwining 
operators satisfy differential equations of regular singular points 
\cite{H2}. He also pointed out that 
if a vertex operator algebra $V$ satisfies $C_2$-cofiniteness 
($\dim V/C_2(V)<\infty$) and rationality (all modules are 
completely reducible), then the associativity of 
product of intertwining operators, or a weaker version which, 
in physicists' terminology, is called the (nonmeromorphic) 
operator product expansion of chiral vertex operators, and 
consequently the category of $V$-modules has a natural structure of 
vertex tensor category and braided tensor category. Namely, their 
results give an extension of tensor product theory 
(cf. \cite{KL}, \cite{TK}) derived from the Knizhnik-Zamolodchikov 
equation \cite{KZ} for the Wess-Zumino-Novikov-Witten models 
\cite{W} and the Belavin-Polyakov-Zamolodchikov equation \cite{BPZ} 
for the minimal models. 

However, in order to apply the tensor product 
theory developed by Huang and Lepowsky, we need assumptions of 
compatibility condition and local 
grading restriction (equivalent to $L(0)$-semisimplicity) 
or rationality (completely reducibility of all modules). Regretfully, 
different from the study of objective models, 
in the axiomatical research for a vertex operator algebra, 
the proof of completely reducibility for all modules is one of 
the hard problems. In the study of any algebras, 
the importance of tensor product theory as a tool is disputable. 
If the category of modules has a natural tensor structure, a tensor 
product becomes a powerful tool and has immeasurable applications. 
Therefore, we need the tensor product theory of vertex operator 
algebra without assuming the rationality. However, there is one 
problem. When we don't assume the rationality, different from the 
classical models in physics, we can't derive the semisimplicity of 
grading operator $L(0)$ from the axioms. Actually, in the recent 
paper in physics, they constructed a few examples, in which $L(0)$ 
does not act semisimply on a module.  For example, it is easy to 
check that a triplet model with central charge $-2$ satisfies 
$C_2$-cofiniteness 
(and so it has only finitely many irreducible modules), but some 
module has an element $v$ satisfying $L(0)v\not=0$ and $L(0)^2v=0$ 
(see \cite{GK}). Furthermore, the differential equations derived from 
the classical differential equations like KZ-equations and 
BPZ-equations, and also derived from $C_1$-cofiniteness, are all of 
regular singular points and so the general solutions have logarithmic 
forms. This fact suggests it natural to include logarithmic 
intertwining operators (and so modules on which $L(0)$ does not act 
semisimply) as the author did in \cite{Miy}. 
The author also showed that $C_2$-cofiniteness is 
more essential than the rationality in \cite{Miy} when 
we consider the models with finitely many simple modules. 
Based on these results, there are several papers from a view point 
of $C_2$-cofiniteness (eg \cite{GK},\cite{NT}). Moreover, the 
confirmation of $C_2$-cofiniteness is discernibly easier than that 
of rationality. For example, if $V$ contains a subVOA $W$ with the 
same Virasoro element and $W$ is $C_2$-cofinite, then so is $V$. 

The aim of the first part of this paper is, following to the paper 
about the modular invariance on the trace function \cite{Miy}, 
to demonstrate how the theory of tensor products develops naturally 
if we include 
{\bf Artin module} (a weak module with a composition series of 
finite length) and {\bf logarithmic intertwining operators}. 
We believe that they are adequate to be called 
modules and intertwining operators of vertex operator algebras. \\

\noindent  
{\bf Main Theorem}[Th 4.5, Th 5.4]\qquad {\it If $V$ is 
$C_2$-cofinite, then for any Artin $V$-modules $W_1$ and $W_2$, 
there is a tensor product module $W_1\boxtimes_{\la}W_2$ for each 
$0\not=\la\in \C$ such that $W_1\boxtimes_{\la}W_2$ is an Artin 
$V$-module and is uniquely determined up to isomorphism. 
Furthermore, the tensor products have natural properties of 
associativity and commutativity. }\\

We should note that let $J_{L(0)}(W_1\boxtimes_{\la}W_2)$ be the 
$L(0)$-radical (see Def. 2.5), then 
$W_1\boxtimes_{\la}W_2/J_{L(0)}(W_1\boxtimes_{\la}W_2)$ is a 
$P(\la)$-tensor product in a sense of \cite{HL1}, that is, it is 
an ordinary $V$-module satisfying the universal property of 
tensor product. \\

In the second part of the paper, we will show an application of our 
theory 
to emphasize the advantages of $C_2$-cofiniteness. For example, 
there is a conjecture: "If $g$ is an automorphism of $V$ of finite 
order, then every irreducible $V^{\langle g\rangle }$-module is 
contained in a $h$-twisted $V$-module for some 
$h\in \langle g\rangle $". We will prove this 
if $V^{\langle g\rangle }$ is $C_2$-cofinite. \\

\noindent
{\bf Theorem 6.9} \quad 
{\it Let $V$ be a simple vertex operator algebra and $g$ an 
automorphism of $V$ of finite order. If a subVOA 
$V^{\langle g\rangle }$ of fixed points is $C_2$-cofinite, then 
every irreducible $V^{\langle g\rangle }$-module is contained in a 
$h$-twisted $V$-module for some $h\in \langle g\rangle $.} \\

\noindent
{\bf Theorem 6.11} \quad 
{\it Let $V$ be a simple vertex operator algebra and 
$G$ an automorphism group of $V$ of finite order. If a subVOA $V^G$ 
of fixed points is rational and $C_2$-cofinite, then $V$ is 
$g$-rational and $V^{\langle g\rangle }$ is rational for every 
$g\in G$. }\\

The many parts in this paper are just to make sure that the theory 
developed by Huang and Lepowsky is still valid without assuming 
the compatiblity condition nor the semisimplicity of grading 
operator. The most part follow from the existence of 
differential equations. For example, 
they introduced a $P(z)$-intertwining 
map for an intertwining operator $\CY(v,z)$ and considered the weight decomposition. 
However, since this is just a value of an intertwining operator $\CY$ at $z$, 
it doesn't concern if $\CY$ is logarithmic or not as long as it exists as a solution 
of differential equations. Moreover, 
$C_2$-cofiniteness guarantees the existence 
of weight space (as a generalized eigenspace of $L(0)$).    

The present paper is organized as follows: In \S 1-\S5, we develop 
the tensor product theory under the assumption of $C_2$-cofiniteness 
on $V$. In \S 6, we introduce a concept of projective modules and 
show rationality of some orbifold models (fixed pointed 
subVOAs by finite automorphisms) as an application of our theory.

\section{Fundamental concepts and their extensions}
In this paper, for a vector space $W$, $W^{\ast}$ denotes the dual 
space $\Hom(W,\C)$ and for a $\N$-graded vector space 
$W=\oplus_{n\in \N}W(n)$, $W'$ denotes the restricted dual space \\
$\oplus_{n\in \Z}(\Hom(W(n),\C))$ and $\ol{W}=\prod_{n\in \N}W(n)$. 

\begin{dfn}[\cite{B}]
A {\bf {\it vertex algebra}} (shortly VA) is a $\Z$-graded vector 
space (graded by weights) $V=\coprod_{n\in \Z}V_{(n)} $
equipped with a linear map 
$$\begin{array}{rl}
V \to &\End(V)[[z,z^{-1}]] \cr
v \to &Y(v,z)=\sum_{n\in \Z}v(n)z^{-n-1}  \qquad
( \mbox{where  }v(n)\in \End(V)), \end{array}$$
$Y(v,z)$ is called the vertex operator associated with $v$, and 
equipped also with a special vector ${\bf 1}\in V_{(0)}$ 
(the vacuum). 

The following conditions are assumed for $u,v\in V$: \\
the lower truncation condition:   $v(n)u=0$  for $n$ sufficiently 
large, \\
vacuum condition:   $Y(\1,z)=1$ (that is, the identity operator);\\
the creation property:  $Y(v,z)\1\in V[[z]]$ and 
$\lim_{z\to 0} Y(v,z)\1=v$,\\ 
the Jacobi identity (the main axiom): 
$$\begin{array}{l}
z_0^{-1}\delta\left(\frac{z_1-z_2}{z_0}\right)Y(u,z_1)Y(v,z_2)
-z_0^{-1}\delta\left(\frac{z_2-z_1}{-z_0}\right)
Y(v,z_2)Y(u,z_1)\cr
=z_2^{-1}\delta\left(\frac{z_1-z_0}{z_2}\right)Y(Y(u,x_0)v,x_2)
\end{array}$$
\end{dfn}

We will consider a vertex algebra satisfying restricted conditions.

\begin{dfn}[\cite{FLM}]
A {\bf {\it vertex operator algebra}} (shortly VOA) is a vertex 
algebra 
$$ V=\oplus_{n=-k}^{\infty}V_{(n)}, \qquad (\dim V_{(n)}<\infty)  $$
and it has a special element $\omega\in V_{(2)}$ 
(the Virasoro element) satisfying \\
the Virasoro algebra relations : 
$[L(m),L(n)]=(m-n)L(m+n)+\frac{1}{12}(m^3-m)\delta_{n+m,0}c$,\\  
where $L(n)$ denotes $\omega(n+1)$ and $c\in \C$ is called central 
charge of $V$,\\
the grading operator:    $L(0)|_{V_{(n)}}=n$ for all $n\in \N$,\\
the $L(-1)$-derivative property: $Y(L(-1)v,z)=\frac{d}{dz}Y(v,z)$ 
for all $v\in V$. 
\end{dfn}

Throughout this paper, $V$ always denotes a vertex 
operator algebra $(V,Y,\1,\omega)$. \\

The following properties are fundamental. 
For $v$, $v_1$, $v_2\in V$ and $v'\in V'$, there is a unique 
$f\in \C[x_1^{\pm},x_2^{\pm},(x_1-x_2)^{-1}]$ such that 
$$\begin{array}{l}
\langle v', Y(v_1,x_1)Y(v_2,x_2)v\rangle 
=\iota_{x_1,x_2}f(x_1,x_2)  \cr
\langle v', Y(v_2,x_2)Y(v_1,x_1)v\rangle 
=\iota_{x_2,x_1}f(x_1,x_2)  \qquad \mbox{  and}\cr
\langle v', Y(Y(v_1,x_1-x_2)v_2,x_2)v\rangle 
=\iota_{x_2,x_1-x_2}f(x_1,x_2).
\end{array}$$
Here, $\iota_{x_1,x_2}f(x_1,x_2)$ denotes an expansion of 
$f(x_1,x_2)$ into a formal series with powers of $x_1$ and $x_2$ 
such that it contains at most finitely negative powers of $x_2$. 
(cf. Proposition 2.3 and 2.4 in \cite{HL1}.)

\begin{dfn}
A {\bf weak module} for VOA $(V, Y, {\bf 1}, \omega)$ is a 
vector space $M$; equipped with a formal power series 
$$Y^M(v, z)=\sum_{n \in {\bf Z}} v^M(n)z^{-n-1} 
\in ({\rm End}(M))[[z, z^{-1}]] $$ 
(called {\it the module vertex operator} of $v$) for $v \in V$ 
satisfying: \\
{\rm (1)}\quad  $v^M(n)w=0$ for $n>\!\!>0$ 
where $v\in V$ and $w\in M$; \\
{\rm (2)}\quad  $Y^M({\bf 1}, z)=1_M$; \\
{\rm (3)}\quad  $Y^M(\omega, z)
=\sum_{n\in\Z} L^M(n)z^{-n-2}$ satisfies: \\
$ \mbox{\ }$ \quad{\rm (3.a)}  the Virasoro algebra relations: 
$$ [L^M(n), L^M(m)]
=(n\!-\!m)L^M(n\!+\!m)+ \delta_{n+m, 0}\frac{n^3-n}{12}c, 
\quad c\in \C $$
$ \mbox{\ }$ \quad{\rm (3.b)}  the $L(-1)$-derivative property:
$  Y^M(L(-1)v, z)=\frac{d}{dz}Y^M(v, z), $ \\
{\rm (4)}\quad  "Jacobi identity" holds: for any $v, u\in V$ and 
$w\in M$, we have 
$$\begin{array}{l}
z_0^{-1}\delta\left(\frac{z_1-z_2}{z_0}\right)Y^M(u,z_1)Y^M(v,z_2)w
-z_0^{-1}\delta\left(\frac{z_2-z_1}{-z_0}\right)
Y^M(v,z_2)Y^M(u,z_1)w\cr
=z_2^{-1}\delta\left(\frac{z_1-z_0}{z_2}\right)Y^M(Y(u,x_0)v,x_2)w
\end{array}$$
\end{dfn}

An {\bf ordinary module} is a $\Z_+$-graded weak $V$-module 
$W=\oplus_{m=0}^{\infty}W(m)$ satisfying \\
{\rm (1')}\quad  if $v\in V_{(r)}$ then $v^W(m)W(n)\subseteq 
W(n\!+\!r\!-\!m\!-\!1)$ \qquad and \\
{\rm (2')}\quad  $\dim W(n)<0$ \qquad and \\
{\rm (3')}\quad  $L(0)$ acts on $W(n)$ semisimply. 

For a simple ordinary module $W=\oplus_{i=0}^{\infty}W(i)$, 
the eigenvalue of $L(0)$ on $W(0)$ is called a 
{\it conformal weight} of $W$. \\

For a weak module $W$, ${\rm soc}_1(W)$ denotes the direct sum of 
ordinary simple modules and we define ${\rm soc}_n(W)$ by 
${\rm soc}_n(W)/{\rm soc}_{n-1}(W)={\soc}_1(W/{\rm soc}_{n-1}(W))$ 
inductively. 
If ${\rm soc}_n(W)=W$, we call {\bf Loewy length} of $W$ is $n$.

One of the main central characters in this paper is 
the following module. 

\begin{dfn}
An {\bf Artin module} is a weak $V$-module W which carries a 
composition series of finite length such that each composition 
factor is a simple ordinary module. 
\end{dfn}

Therefore, an Artin module $W$ has also a $\Z_+$-grading 
$W=\oplus_{n=0}^{\infty}W(n)$ satisfying (1'). If $w\in W(n)$, 
we denote it by $\deg(w)=n$. Since $L(0)$ acts on $W(n)$ for 
each $n$, $W$ is also a direct sum of finite dimensional generalized 
eigenspaces $W_{(r)}=\{w\in W\mid (L(0)-r)^nw=0 
\mbox{  for some }n\in \N\}$. We call a nonzero element $w$ in 
$M_{(r)}$ {\rm homogeneous element} with weight $r$ and denote it 
by $\wt(w)=r$. For $v\in W_{(r)}$, we will call the minimal integer 
$n$ satisfying $(L(0)-r)^nv=0$ 
{\bf the $L(0)$-Loewy length} of $v$ by $L(0)$.

From the axioms, the grading operator $L(0)$ commutes with all 
grade-preserving operators. Therefore we can define: 

\begin{dfn}
For an Artin module $W=\oplus_{r\in \C}W_{(r)}$, 
define 
$$J_{L(0)}(W)=\oplus_{r\in \C}(L(0)-r)W_{(r)},$$ 
then $U$ is an Artin submodule and $W/J_{L(0)}(W)$ is the largest 
$L(0)$-semisimple factor module of $V$. 
We will call $J_{L(0)}(W)$ a {\bf $L(0)$-radical} of $W$. 
\end{dfn}

Clearly we can define a restricted dual module for an Artin module. 

\begin{dfn}
Let $(W,Y^W)$ be an Artin module. A vertex operator $Y'(v,z)$ on 
$W'=\coprod_{n\in \N}\Hom(W(n),\C)$ is defined by 
$$\langle Y'(v,z)w',w\rangle =\langle w', 
Y^W(e^{L(1)z}(-z^{-2})^{\wt(v)}v, z^{-1})w\rangle  \eqno{(2.1)}$$
for $w'\in W', w\in W, v\in V_{(\wt(v))}$ and extend it linearly. 
Then $(W',Y')$ is an Artin $V$-module and we call it 
{\bf a restricted dual} of $W$.
\end{dfn}

We denote $Y^W(e^{L(1)z}(-z^{-2})^{\wt(v)}v, z^{-1})$ by 
$Y^{\ast}(v,z)$. 

\begin{dfn} 
For a $V$-module $W$, set 
$$\begin{array}{l}
C_1(W)=\langle v(-1)u \mid v\in V, u\in W, \wt(v)\geq 1\rangle, 
\cr
C_2(W)=\langle v(-2)u \mid v\in V, u\in W, \wt(v)\geq 0\rangle. 
\end{array} \eqno{(2.2)}$$ 
We call $W$ $C_n$-cofinite if $\dim W/C_n(W)<\infty$. 
\end{dfn}

If $W$ is 
$C_n$-cofinite, then we have a finite subset $B_n(W)$ of $W$ of 
homogeneous elements such that $W=C_n(W)+\langle B_n(W)\rangle$. 

\begin{rmk}
If $V$ is $C_2$-cofinite, then the following are known: \\
(i) every weak module is $\Z_+$-graded and every finitely generated 
weak module is an Artin module, {\rm \cite{Miy}}.\\
(ii) $n$-th Zhu algebra $A_{1,n}(V)$ is finite dimensional and the 
number of inequivalent irreducible modules is finite, 
{\rm \cite{GN},\cite{DLM}}.\\
(iii) conformal weights are all rational numbers, {\rm \cite{Miy}}.\\
(iv) every Artin module is $C_2$-cofinite, 
{\rm \cite{Bu}, \cite{GN}}.\\
(v) for any weak module W generated from one element $w$ has the 
following spanning set $\{v_1(n_1)....v_k(n_k)w \mid v_i\in B_2(V), 
\quad  n_1<\cdots <n_k\}$, {\rm \cite{Miy}}.   
\end{rmk}

Although the ordinary definition of $C_2(V)$ introduced by Zhu 
\cite{Z} is 
$$C_2(W)=\langle v(-2)u \mid v\in V, u\in W\rangle \eqno{(2.3)} $$
and the above results were proved under this definition, 
we adopt (2.2) in order to cover the case 
where $V$ has negative weights. For, in their proof, 
the advantage of $C_2$-cofiniteness is just to satisfy 
$\wt(u(-2)v)\gneq \wt(u),\wt(v),\wt(u)+\wt(v)$ so that they could 
use an induction.

\subsection{Intertwining operators for logarithmic modules}
In \cite{Mil}, A.~Milas extended a concept of intertwining operators 
defined by \cite{FHL}. Namely, he introduced a logarithmic 
intertwining operators, which is one of the main characters 
in this paper.

\begin{dfn} Let $W_1$, $W_2$ and $W_3$ be Artin $V$-modules. 
A (logarithmic) intertwining operator of type 
$\binom{W_3}{W_1W_2}$ is a linear map 
$$\begin{array}{l}
  \CY(,x):W_1 \to \Hom(W_2,W_3)\{z\}[\log(z)] \cr
  \CY(w,x)=\sum_{i=0}^K\sum_{n\in \C} w_{(n,i)}z^{-n-1}(\log(z))^i 
\end{array}$$
satisfying the following conditions: \\
1. The lower truncation property: 
(for each coefficient $\sum_{n\in C}w_{(n,i)}z^{-n-1}$ of 
$(\log(z))^i$) \\ 
2. $L(-1)$-derivative condition: 
$\CY(L(-1)w,z)=\frac{d}{dz}\CY(w,z)$ \\
3. The Jacobi identity:   In $W_3\{x_0,x_1,x_2\}[\log(x_2)]$, 
for $u\in V$, $w_1\in W_1$, $w_2\in W_2$,  
$$\begin{array}{l}
x_0^{-1}\delta(\frac{x_1-x_2}{x_0})Y(u,x_1)\CY(w_1,x_2)w_2
-x_0^{-1}\delta(\frac{x_2-x_1}{-x_0})Y(w_1,x_2)\CY(u,x_1)w_2 \cr
=x_2^{-1}\delta(\frac{x_1-x_0}{x_2})\CY(Y(u,x_0)w_1,x_2)w_2 
\end{array}$$
We denote the space of logarithmic intertwining operators of type 
$\binom{W_3}{W_1W_2}$ by $\CV^{W_3}_{W_1W_2}$. 
\end{dfn}

A logarithmic intertwining operator satisfies, just like an ordinary 
one, \\
(1) {\bf Commutativity}: $[v(n), \CY(w,z)]=
\displaystyle{\sum_{i=0}^{\infty}} \binom{n}{i}\CY(v(i)w,z)z^{n-i} 
\mbox{  for  }v\in V$ and\\
(2) {\bf Associativity}: $\CY(v(-1)w,z)=
\displaystyle{\sum_{m=0}^{\infty}}v(-m-1)z^m\CY(w,z)
+\CY(w,z)\displaystyle{\sum_{m=0}^{\infty}}v(m)z^{-m-1}$ \\
for $v\in V$. For a logarithmic intertwining operator 
$\CY(w,z)=\displaystyle{\sum_{k=0}^K\sum_{n\in \C}}
w_{(n,k)}z^{-n-1}(\log(z))^k$, we denote the coefficient 
$\sum_{n\in \C}w_{(n,k)}z^{-n-1}$ of $(\log(z))^k$ by 
$\tilde{\CY}^{(k)}(w,z)$. If $\tilde{\CY}^{(K)}(w,z)\not=0$ for some 
$w\in W_1$, then we call that {\bf $L(0)$-Loewy length} of $\CY$ is 
$K+1$. Since $Y(u,x)$ does not contain $\log(x)$-terms, each 
$\tilde{\CY}^{(i)}(w,x)$ also satisfies all conditions except for 
$L(-1)$-derivative condition. On the other hand, 
from $L(-1)$-derivative condition, we have 
$$ \begin{array}{l}
[L(-1),w_{(n,i)}]=-nw_{(n-1,i)}+(i+1)w_{(n-1,i+1)} \cr
L(0)w_{(n,i)}u=(L(0)w)_{(n,i)}u+w_{(n,i)}(L(0)u)
+(-n-1)u_{(n,i)}+(i+1)w_{(n,i+1)}u. 
\end{array} \eqno{(2.4)}$$
If $w$ and $u$ are homogeneous, then $w_{(n,i)}u$ is also 
a homogeneous element and 
$$ \wt(w_{(n,m)}u)=\wt(w)+\wt(u)-n-1. $$
Furthermore, if $w$ and $u$ are eigenvectors of $L(0)$, then 
$$L(0)w_{(n,i)}u=(\wt(w)+\wt(u)-n-1)w_{(n,i)}u+(i+1)w_{(n,i+1)}u.$$

\vspace{2mm}

\begin{rmk} If $W_1$,$W_2$ and $W_3$ are all ordinary modules, 
then any logarithmic intertwining operator of type 
$\binom{W_3}{W_1W_2}$ is an ordinary intertwining operator. 
Therefore, the fusion rules $\dim \CV^{W_3}_{W_1W_2}$ among 
three irreducible modules agree with the classical fusion rules. 
\end{rmk}

We also note that powers of $z$ in 
logarithmic intertwining operators are rational numbers if $V$ 
is $C_2$-cofinite since 
all conformal weights are rational numbers. For a complex 
number $\la$ and a real number $n$, ${\la}^n$ always 
denotes $e^{n\log(\la)}$ in this paper, where 
$\log(\la)=\log|\la|+i \arg(\la)$ with $0\leq \arg(\la)<2\pi$.

\section{$n$-point functions and differential equations}
The purpose of this section is just to confirm that the differential 
equations which was introduced by Huang \cite{H2} and satisfied by 
$n$-point correlation functions of intertwining operators under 
the condition that all modules are $C_1$-cofinite are still valid 
for logarithmic intertwining operators.

Let $V=\coprod_{n\in \Z}V_{(n)}$ be a vertex operator algebra, 
$R=\C[z_1^{\pm 1}, z_2^{\pm 1}, (z_1-z_2)^{-1}]$ and define a 
filtration $R_{-m}=(z_1)^{-m}\C[z_1,z_2]+z_2^{-m}\C[z_1,z_2]+
(z_1-z_2)^{-m}\C[z_1,z_2]  \quad (m=0,1,...)$.

Let $W_0$,$W_1$,$W_2$ and $W_3$ be Artin modules satisfying 
$C_1$-cofiniteness. For $i=0,1,2,3$, let $p_i$ be natural 
numbers satisfying $W_i=C_1(W_i)+\sum_{m=0}^{p_i} W_i(m)$ and 
set $p=\sum_{i=0}^3p_i$. We use the notation $\deg(w)=m$ to denote 
$w\in W_i(m)$ and extend it to the degree on 
$W_0\otimes ...\otimes W_3$ by total degree. Set 
$$W(m)=\sum_{m=r_0+r_1+r_2+r_3}(W_0(r_0)\otimes 
W_1(r_1)\otimes W_2(r_2)\otimes W_3(r_3)).$$

Let $W_4$ and $W_5$ be also Artin $V$-modules. For logarithmic 
intertwining operators $\CY_1$, $\CY_2$, $\CY_3$ and $\CY_4$ 
of type $\binom{W_0'}{W_1W_4}$, $\binom{W_4}{W_2W_3}$, 
$\binom{W_5}{W_1W_2}$ and $\binom{W_0'}{W_5W_3}$, respectively, 
define $2$-point correlation functions (formal series) by 
$$\begin{array}{l}
S^1(w_0\otimes w_1\otimes w_2\otimes w_3;z_1,z_2)
=\langle w_0, \CY_1(w_1,z_1)\CY_2(w_2,z_2)w_3\rangle,  \cr
S^2(w_0\otimes w_1\otimes w_2\otimes w_3;z_1,z_2)=
\langle w_0, \CY_4(\CY_3(w_1,z_1-z_2)w_2,z_2)w_3\rangle\end{array}$$ 
and for $f(z_1,z_2)\in R, u\in W_0\otimes ...\otimes W_4$, 
extend it $R$-linearly, that is,   
$S^i(u\otimes f(z_1,z_2))=S^i(u)f(z_1,z_2)$ for $f(z_1,z_2)\in R$. 
In the following arguments, since we can treat $S^2$ as well as 
$S^1$, we will prove the assertion only for $S^1$. We also define 
$$S^{1,(i,j)}(w_0\otimes w_1\otimes w_2\otimes w_3;z_1,z_2)=
\langle w_0, \tilde{\CY}_1^{(i)}(w_1,z_1)\tilde{\CY}_2^{(j)}(w_2,z_2)w_3\rangle $$
and $S^{1-}$ denotes one of $S^1, S^{1,(i,j)}$. Since every 
$\tilde{\CY}_i^{(m)}$ satisfies Commutativity and Associativity with 
vertex operators on modules, direct calculations show:
$$ \begin{array}{l}
S^{1-}(v^{\ast}(m)w_0\otimes w_1\otimes w_2\otimes w_3;z_1,z_2)\cr
\mbox{}\quad =\sum_{k=0}^{\infty}\binom{m}{k}(z_1)^{m-k}
S^{1-}(w_0\otimes v(k)w_1\otimes w_2\otimes w_3;z_1,z_2) \cr
\mbox{}\qquad +\sum_{k=0}^{\infty}\binom{m}{k}(z_2)^{m-k}
S^{1-}(w_0\otimes w_1\otimes v(k)w_2\otimes w_3;z_1,z_2) \cr 
\mbox{}\qquad +
S^{1-}(w_0\otimes w_1\otimes w_2\otimes v(m)w_3;z_1,z_2), \cr
S^{1-}(w_0\otimes v(-1)w_1\otimes w_2\otimes w_3;z_1,z_2)\cr
\mbox{}\quad =\sum_{k=0}^{\infty}(z_2)^k
S^{1-}(v^{\ast}(-1-k)w_0\otimes w_1\otimes w_2\otimes w_3;z_1,z_2)\cr
\mbox{}\qquad +\sum_{k=0}^{\infty}\binom{-1}{k}(z_1+z_2)^{-1-k}
S^{1-}(w_0\otimes w_1\otimes v(k)w_2\otimes w_3;z_1,z_2) \cr 
\mbox{}\qquad +\sum_{k=0}^{\infty}(z_1)^{-1-k}
S^{1-}(w_0\otimes w_1\otimes w_2\otimes v(k)w_3;z_1,z_2), \cr
S^{1-}(w_0\otimes w_1\otimes v(-1)w_2\otimes w_3;z_1,z_2)\cr
\mbox{}\quad =\sum_{k=0}^{\infty}(z_2)^k
S^{1-}(v^{\ast}(-1-k)w_0\otimes w_1\otimes w_2\otimes w_3;z_1,z_2)\cr
\mbox{}\qquad -\sum_{k=0}^{\infty}\binom{-1}{k}(z_1-z_2)^{-1-k}
S^{1-}(w_0\otimes v(k)w_1\otimes w_2\otimes w_3;z_1,z_2) \cr 
\mbox{}\qquad +\sum_{k=0}^{\infty}(z_2)^{-1-k}
S^{1-}(w_0\otimes w_1\otimes w_2\otimes v(k)w_3;z_1,z_2), \cr
S^{1-}(w_0\otimes w_1\otimes w_2\otimes v(-1)w_3;z_1,z_2)\cr
\mbox{}\quad =
S^{1-}(v^{\ast}(-1)w_0\otimes w_1\otimes w_2\otimes w_3;z_1,z_2) \cr
\mbox{}\qquad -\sum_{k=0}^{\infty}\binom{-1}{k}(z_1)^{-1-k}
S^{1-}(w_0\otimes v(k)w_1\otimes w_2\otimes w_3;z_1,z_2) \cr 
\mbox{}\qquad -\sum_{k=0}^{\infty}(z_2)^{-1-k}
S^{1-}(w_0\otimes w_1\otimes v(k)w_2\otimes w_3;z_1,z_2)  
\end{array}\eqno{(3.1)}$$
for $v\in V, w_0\in W_0, w_1\in W_1, w_2\in W_2, w_3\in W_3$, 
where $v^*(m)$ denotes the dual operator of $v(m)$. 
We note that every formal infinite sum becomes a finite sum, 
since $v(k)w_{i}=0=v^*(-k)w_0$ for $k$ sufficiently large.

\begin{lmm}
Under the above setting, for any $w \in W(m)$ and $i=1,2$, 
$S^i(w;z_1,z_2)$ is expressed by a sum of functions of the form 
$$ S^i(u^{(k)};z_1,z_2)f_k \quad u^{(k)}\in W(k), k\leq p, f_k 
\in R_{-k+m}. $$ 
\end{lmm}

\pr 
We may assume that $w$ has a form 
$u_0\otimes u_1\otimes u_2\otimes u_3$ and $u_i$ are homogeneous. 
We will prove the lemma by the induction on 
$m=\sum_{i=0}^3 \deg(u_{i})$. 
If $\sum_{i=0}^3 \deg(u_i)\leq p$ then the assertion is trivial. 
If $\sum_{i=0}^3 \deg(u_i)>p$ then some $u_i$ is in $C_1(W_i)$ by 
the definition of $p$. Therefore, there are homogeneous elements 
$v_j\in V$ $(\wt(v_j)\geq 0)$ and $u_i^j\in W_i$ such that 
$u_i$ is a sum of elements of the form $v_j(-1)u_i^j$. Therefore, 
we may assume $u_i=v(-1)w_i$. In (3.1), rewriting $u_j$ by $w_j$ 
for $i\not=j$, we have shown that 
$S^k(w_0\otimes \cdots \otimes v(-1)w_i\otimes \cdots \otimes w_4)$
is a linear sum of 
$$ \begin{array}{ll}
S^k(v(-n)^{\ast}w_0\otimes w_1\otimes \cdots \otimes w_i\otimes 
\cdots \otimes w_4)f_1, \quad &n\geq 1, f_1\in R_{2a+n-1}\cr
S^k(w_0\otimes \cdots v(n)w_2 \otimes \cdots w_i\otimes \cdots 
\otimes w_4)f_2, \quad &n\geq 0, f_2\in R_{n+1}\cr
S^i(w_0\otimes \cdots \otimes v(n)u\otimes \cdots)f_3, 
\quad &n\geq 0, f_3\in R_{n+1}. 
\end{array}$$ 
Since $\wt(v(n)), \wt(v^{\ast}(-n)) \lneq \wt(v(-1))$ for $n\geq 0$ 
as operators, the total degree of elements in 
$W_0\otimes \cdots \otimes W_3$ which appear in the right-hand side 
of (3.1) is less than $\deg(v)+\sum_{i=0}^3 \deg(w_i)$. Therefore, 
by using the induction, we have the desired result. 
\prend

Using the above result, we can prove the existence of differential 
equations of regular singular points for $2$-point correlation 
functions. 

\begin{thm}[cf. Th.1.4, Th.2.3  in \cite{H2}]
Let $W_i$ for $i=0,1,2,3$ be Artin modules satisfying 
$C_1$-cofiniteness. Then for any $w_i\in W_i (i=0,1,2,3)$, there 
exist 
$$  a_k(z_1,z_2), b_l(z_1,z_2) \in \C[z_1^{\pm 1}, z_2^{\pm 1}, 
(z_1-z_2)^{-1}]$$
for $k=1,2,...,m$ and $l=1,...,n$ such that for any $V$-modules 
$W_4$, $W_5$ and $W_6$, any (logarithmic) intertwining operators 
$\CY_1,...,\CY_6$ of types $\binom{W_0'}{W_1W_4}$, 
$\binom{W_4}{W_2W_3}$, $\binom{W_5}{W_1W_2}$, 
$\binom{W_0'}{W_5W_3}$, $\binom{W_0'}{W_2W_6}$ and 
$\binom{W_6}{W_1W_3}$, respectively, the series 
$$\langle w_0,\CY_1(w_1,z_1)\CY_2(w_2,z_2)w_3\rangle\eqno{(3.2)}$$
$$\langle w_0,\CY_4(\CY_3(w_1,z_1-z_2)w_2,z_2)w_3\rangle  
\eqno{(3.3)}$$
$$\langle w_0,\CY_5(w_2,z_2)\CY_6(w_1,z_1)w_3\rangle \eqno{(3.4)}$$
satisfy the expansion of the system of differential equations of 
regular singular points 
$$   \frac{\pd^m\varphi}{\pd z_1^m}+a_1(z_1,z_2)
\frac{\pd^{m-1}\varphi}{\pd z_1^{m-1}}+\cdots+a_m(z_1,z_2)
\varphi=0 \eqno{(3.5)}$$
$$   \frac{\pd^n\varphi}{\pd z_2^n}+b_1(z_1,z_2)
\frac{\pd^{n-1}\varphi}{\pd z_1^{n-1}}+\cdots+b_n(z_1,z_2)
\varphi=0 \eqno{(3.6)}$$
in the regions $|z_1|>|z_2|>0$, $|z_2|>|z_1-z_2|>0$ and 
$|z_2|>|z_1|>0$, respectively. 
\end{thm}

\pr
We will restrict ourself to the case of (3.2), 
since we can prove the assertions for (3.3) and (3.4) similarly. 
Set $P=\oplus_{i=0}^p W(m)$. Then $P\otimes R$ is a finitely 
generated $R$-module and the derivation 
$\frac{d}{dz_i}$ acts on this space as an $R$-linear operator 
$$ \begin{array}{l}
 \frac{d}{dz_1}S^1(w_0\otimes w_1\otimes w_2 \otimes w_3;z_1,z_2)
=S^1(w_0\otimes L(-1)w_1\otimes w_2\otimes w_3;z_1,z_2)\cr
 \frac{d}{dz_2}S^1(w_0\otimes w_1\otimes w_2 \otimes w_3;z_1,z_2)
=S^1(w_0\otimes w_1\otimes L(-1)w_2\otimes w_3;z_1,z_2). 
\end{array} $$

For $k=1,2$, decompose each matrix representation $A^{(k)}$ 
of $\frac{d}{dz_k}$ on $R\otimes P$ into 
$A^{(k)}=(A^{(k)}_{ij})_{i,j=0,...,p}$ 
with small matrices $A^{(k)}_{ij}$ $(i,j=0,...,p)$, 
where $A_{ij}$ is a $\dim W(i)$ by $\dim W(j)$-matrix to denote 
a composition map of $\frac{d}{dz_k}: W(i) \to W$ and the 
restriction $W \to W(j)$. 

Then since the operator $L(-1)$ upgrades the degree just by one, 
the above lemma suggests that we are able to choose  
$$\begin{array}{l}
  A^{(k)}_{i,i+1}\in M(R_0)  \cr
  A^{(k)}_{n,1}\in M(R_{-n}) \cr  
  A^{(k)}_{i,j}=0 \quad \mbox{for else,}
\end{array}$$ 
where $M(R_m)$ denotes the ring of all matrices over $R_m$ of 
given size. Then it is easy to see that the characteristic 
polynomial of $A^k$ (and also of $\frac{d}{dz_k}$) 
$$f^{(k)}(x)=|xI-A^{(k)}|=x^t+\sum_{i=1}^{t} g_i^{(k)}
(z_1,z_2)x^{t-i}$$ satisfies $g_i^{(k)}(z_1,z_2)\in R_{-i}$. 
Therefore, a differential equation $f^{(k)}(\frac{d}{dz_k})=0$ has 
singular points at $z_1=0$, $z_2=0$ and $z_1=z_2$ at most and 
is of regular singular points. 

This completes the proof of Theorem.
\prend

\begin{rmk}
Although we have shown that $2$-point functions satisfy differential 
equations of regular singular points, 
it is clear the same argument works for any $n$-point functions. 
\end{rmk}

As a simple application, we have:

\begin{thm}[cf. \cite{H2}]  Let $V$ be a vertex operator algebra 
and let $W_1$, $W_2$ and $W_3$ be $C_1$-cofinite Artin modules. 
Then the fusion rule among these three modules is finite. 
\end{thm}

\pr
For any logarithmic intertwining operator \\
$\CY(w,z)=\sum_{i=0}^K\sum_{n\in \C} w_{(n,i)}z^{-n-1}\log(z)^i$, 
$S(w_1,w_2,w_3;z)=\langle w_3, \CY(w_2,z)w_1\rangle $ satisfies 
a differential equation 
$$\{(\frac{d}{dz})^n+\sum_{j=1}^n 
a_j(W_1,W_2,W_3)(\frac{d}{dz})^{n-j})\}S(w_1,w_2,w_3;z)=0 $$
and so the space of solutions is of finite dimension for 
each $(w_1,w_2,w_3)$. On the other hand, the above lemma means that 
$S(w_1,w_2,w_3;z)$ for general triplet $(w_1,w_2,w_3)$ is determined 
by a finite set of $S(u_1,u_2,u_3;z)$ satisfying 
$\sum_{i=1}^3 \deg(u_i)\leq p.$ 
Combining them, we have the desired result. \hfill \qed

\section{$P(\la)$-tensor products and $\la$-tensor products}
In this section, we fix a complex number $\la\not=0$. 
The aim of this section is that if $V$ is $C_2$-cofinite, then 
a tensor product $W_1\boxtimes_{\la} W_2$ always exists as an 
Artin module for any two Artin modules $W_1$ and $W_2$.  
The concepts of $P(\la)$-intertwining map and $P(\la)$-tensor product 
$W_1\boxtimes_{P(\la)}W_2$ of ($L(0)$-semisimple) $V$-modules $W_1$ 
and $W_2$ were introduced in \cite{HL1}.  
As we explained in the introduction, they have treated only 
$L(0)$-semisimple modules and study the existence of 
($L(0)$-semisimple) tensor product $W_1\boxtimes_{P(\la)}W_2$ for 
ordinary modules $W_1$ and $W_2$. In order to carry out 
$L(0)$-semisimple tensor products, they used two conditions,  
$P(\la)$-compatibility and $P(\la)$-grading restriction 
condition. However, the second one is just to check the 
$L(0)$-semisimplicity. The advantage of this paper is, because 
we don't demand $L(0)$-semisimplicity, we don't need to check the 
second condition, which makes our story very easy.  We note that 
even if $W_1$ and $W_2$ are irreducible, 
$W_1\boxtimes_{\la}W_2$ may not be $L(0)$-semisimple. 

For an intertwining operator 
$\CY(w,z)=\sum_{i=0}^K\sum_{n\in C} w_{(n,i)}z^{-n-1}(\log z)^i$ 
of type $\binom{W_3}{W_1W_2}$ and $u\in W_2$, define 
$$F^{\CY}_{\la}(w\otimes u)=
\prod_{n\in \R}(\sum_{i=0}^K \la^{-n-1}(\log(\la))^i(w_{n,i}u) )
\in \ol{W_3}. $$

\begin{dfn} 
For Artin modules $W_1$, $W_2$ and $W_3$ and a complex number 
$\la\not=0$, a {\it $P(\la)$-intertwining map} of type 
$\binom{W_3}{W_1\ W_2}$ is a 
linear map $F:W_1\otimes W_2 \to \ol{W_3}$ such that 
$$\begin{array}{l}
x_0^{-1}\delta(\frac{x_1-\la}{x_0})Y_3(v,x_0)F(w_{(1)}
\otimes w_{(2)})\cr
\mbox{}\qquad=\la^{-1}\delta(\frac{x_1-x_0}{\la})F(Y_1(v,x_0)w_{(1)}
\otimes w_{(2)})+x_0^{-1}\delta(\frac{\la-x_1}{-x_0})F(w_{(1)}
\otimes Y_2(v,x_1)w_{(2)})
\end{array}$$
for $v\in V$, $w_{(1)}\in W_1$, $w_{(2)}\in W_2$. 

We denote the space of $P(\la)$-intertwining maps of type 
$\binom{W_3}{W_1\ W_2}$ by $\CM[P(\la)]^{W_3}_{W_1W_2}$. 
A $P(\la)$-product of $W_1$ and $W_2$ is an Artin module $W_3$ 
equipped with a $P(\la)$-intertwining map $F$ of type 
$\binom{W_3}{W_1\ W_2}$. A morphism from $(W_3,Y_3;F)$ to 
$(W_4,Y_4;G)$ is a module map $\zeta$ from $W_3$ to $W_4$ such 
that $G=\ol{\zeta}\circ F$, where $\ol{\zeta}$ is the natural 
map from $\ol{W_3} \to \ol{W_4}$ uniquely extending $\zeta$. 
\end{dfn}

\begin{dfn}
A $\la$-tensor product of $W_1$ and $W_2$ is a $P(\la)$-product of 
$W_1$ and $W_2$  
$$ (W_1\boxtimes_{\la}W_2, Y_{\la};\boxtimes_{\la}) $$
such that for any $P(\la)$-product $(W_3,Y_3;F)$ of $W_1$ and $W_2$, 
there is a unique morphism $\zeta$ from $W_1\boxtimes_{\la}W_{2}$ to 
$W_3$ such that $\ol{\zeta}\circ \boxtimes_{\la}=F$. 
\end{dfn}

We would like to use the name "$\la$-tensor product", but not 
$P(\la)$-tensor product, to tell the difference.  Later, we will 
construct $\la$-tensor product explicitly for $C_2$-cofinite $V$.
We note that $W_1\boxtimes_{\la}W_2/J_{L(0)}(W_1\boxtimes_{\la}W_2)$ 
is a $P(\la)$-tensor product in a sense of \cite{HL1}.

Let's explain relations between $P(\la)$-intertwining maps and 
logarithmic intertwining operators of the same type. 
If $W_1$,$W_2$ $W_3$ are all $L(0)$-semisimple modules, 
the following is the same as in \cite{HL1}. 

Let $\CY(v,z)=\sum_{i=0}^K \tilde{\CY}^{(i)}(v,z)(\log(z))^i$ be a 
(logarithmic) intertwining operator of type $\binom{W_3}{W_1W_2}$. 
We have a linear map $F_{\CY}^{P(\la)}: W_1\otimes W_2 \to \ol{W_3}$ 
given by 
$$  F_{\CY}^{P(\la)}(w_1\otimes w_2)=\CY(w_1,\la)w_2$$
for all $w_1\in W_1$, $w_2\in W_2$. Since $\CY(v,z)$ satisfies 
Jacobi identity, $F_{\CY}^{P(\la)}$ is a $P(\la)$-intertwining map. 
We note that since we use only Jacobi identity in 
the definition of $P(\la)$-intertwining map, 
$$  F_{\tilde{\CY}^{(m)}}^{P(\la)}(w_1\otimes w_2)=\tilde{\CY}^{(m)}(w_1,\la)w_2$$
is also a $P(\la)$-intertwining map for each $m$.

Conversely, given a $P(\la)$-intertwining map $F$, for $n\in \C$ and 
two homogeneous elements $w\in W_1, u\in W_2$, we define 
$w_{(n,0)}u$ to be the projection of the image of $w\otimes u$ under 
$F$ to the homogeneous subspace $(W_3)_{(\wt(w_1)+\wt(w_2)-n-1)}$ of 
$W_3$ of weight $\wt(w_1)-n-1+\wt(w_2)$ 
multiplied by $\la^{(n+1)}$.  Using this, we define $\tilde{\CY}^{(0)}_F$ 
by $\tilde{\CY}^{(0)}_{F}(w,z)u=\sum_{n\in \C} w_{(n,0)}u z^{-n-1}$ and 
extend it linearly. Then $\tilde{\CY}_F^{(0)}$ satisfies the Jacobi identity. 
If $(w_1)_{(n,0)}w_2$, $w_1$ and $w_2$ are eigenvectors of $L(0)$, 
as they showed in \cite{HL1}, $\tilde{\CY}_F^{(0)}$ also satisfies 
$L(-1)$-derivative property. In the general case, using (2.4), we 
define 
$$w_{(n,1)}u=(L(0)w_{(n,0)}u-(L(0)w)_{(n,0)}u-w_{(n,0)}(L(0)u)
+(n+1)w_{(n,0)}u $$
and set $\tilde{\CY}_F^{(1)}(w,z)=\sum_{n\in \C} w_{(n,1)}z^{-n-1}$. 
Similarly, we define 
$$w_{(n,i+1)}u=\frac{1}{i+1}(L(0)w_{(n,i)}u-(L(0)w)_{(n,i)}u
-w_{(n,i)}(L(0)u)+(n+1)w_{(n,i)}u $$
for $i=1,...$ and set 
$$\tilde{\CY}^{(m)}_{F}(w,z)u=\sum_{n\in \C} w_{(n,m)}uz^{-n-1}.$$ 
Since there is an integer $K$ such that 
$(L(0)-r)^K=0$ on every $L(0)$-eigenspace of 
$W_1\oplus W_2\oplus W_3$ with eigenvalue $r$ for every $r\in \C$, 
the above process ceases with finite steps. Therefore, 
$$\CY_F^{\la}(v,z)=\sum_{m=0}^{\infty} \tilde{\CY}_F^{(m)}(v,z)(\log(z))^m$$
is well-defined and satisfies the Jacobi identity and 
$L(-1)$-derivative property. Thus $\CY_F^{\la}(v,z)$ is an 
intertwining operator. 

Repeating the above steps, if we construct an intertwining map 
$F_{\CY}^{P(\la)}$ from an intertwining operator $\CY$ and we also 
construct intertwining operator $\CY_{F_{\CY}}^{\la}$ from 
$F_{\CY}^{P(\la)}$, then the $L(0)$-Loewy length of the difference 
$\CY-\CY^{\la}_{F_{\CY}^{P(\la)}}$ is less than one of $\CY$. 
Therefore, the above constructions give isomorphisms 
between the space $\CV$ of intertwining operators and the space 
of intertwining maps of the same type. \\

Thus, we have the following: 

\begin{prn}  The correspondence $\CY \to F_{\CY}^{P(\la)}$ is a 
linear isomorphism from the vector space $\CV_{W_1W_2}^{W_3}$ of 
intertwining operators of type $\binom{W_3}{W_1W_2}$ to the vector 
space $\CM[P(\la)]_{W_1W_2}^{W_3}$ of $P(\la)$-intertwining maps of 
type $\binom{W_3}{W_1W_2}$.  \end{prn}

\begin{rmk}
By the above argument, if $W_1$ and $W_2$ are ordinary modules and 
$\CY(v,z)=\sum_{k=0}^K \tilde{\CY}^{(k)}(v,z)(\log(z))^k$ is a logarithmic 
intertwining operator of type $\binom{W_3}{W_1W_2}$, then \\
$\tilde{\CY}^{(0)}(v,z)$ induces an intertwining operator of type 
$\binom{W_3/J_{L(0)}(W_3)}{W_1, W_2}$. In particular, if a 
$\la$-tensor product $W_1\boxtimes_{\la} W_2$ exits, then 
$W_1\boxtimes_{\la} W_2/J_{L(0)}(W_1\boxtimes_{\la}W_2)$ is a 
$P(\la)$-tensor product $W_1\boxtimes_{P(\la)}W_2$ in a sense of 
\cite{HL1}. 
\end{rmk}

For a logarithmic intertwining operator 
$\CY(w,z)=\sum_{k=0}^{K}\sum_{n\in \C}w_{(n,k)}z^{-n-1}(\log(z))^k$ 
of type $\binom{W_3}{W_1W_2}$ and $0\not=\la\in \C$ and 
$u'\in W_3'$, we can define $\phi^{\la}_{\CY}: W_3' \to 
(W_1\otimes W_2)^{\ast}$ by 
$$ \phi^{\la}_{\CY}(u')(w_1\otimes w_2)=\langle u', 
\CY(w_1,\la)w_2\rangle . $$
If $\langle w_{(n,i)}u\mid w\in W_1, u\in W_2, n\in \C, i\in \N
\rangle =W_3$, 
then $\phi^{\la}_{\CY}$ is injective and we can embed $W_3'$ in 
$(W_1\otimes W_2)^{\ast}$. 

Putting the above embedding in our mind, for two Arting $V$-modules 
$(W_1,Y_1)$ and $(W_2,Y_2)$ and $\la\in\C$, we define the actions 
$Y(v,z)=\sum_{n\in \Z} v(n)z^{-n-1}$ of $v\in V$ on $W_1\otimes W_2$ 
by 
$$ v(n)(w_1\otimes w_2)=(\sum_{i=0}^{\infty} \binom{n}{i}
\la^{n-i}v(i)w_1\otimes w_2)+(w_1\otimes (v(n)w_2) \eqno{(4.1)}$$ 
and we denote $W_1\otimes W_2$ with such actions by 
$W_1\otimes_{\la} W_2$. We then define the actions of $V$ on 
$(W_1\otimes_{\la} W_2)^{\ast}$ by 
$\langle Y_{\la}(v,z)f, w_1\otimes_{\la} w_2 \rangle
=\langle f, Y^{\ast}(v,z)(w_1\otimes_{\la} w_2)\rangle $. The above 
formulas do indeed give well-defined maps in generating-function 
forms and the vertex operators on $(W_1\otimes W_2)^{\ast}$ satisfy 
commutativity and $L(-1)$-derivative property if the lower truncation 
property hold, as they showed in \cite{HL2} and \cite{L}. We set  
$$D((W_1\otimes_{\la} W_2)^{\ast})=
\left\{w\in (W_1\otimes_{\la} W_2)^{\ast}\mid 
\begin{array}{l} Y_{\la}(v,z)w\mbox{ is Lorentzian for all } 
v\in V \cr \mbox{Jacobi identity (Def.2.3 (4)) holds on $w$ }
\end{array} \right\},$$
which is a weak module by the definition. \\

We next prove the following:

\begin{lmm}
If $V$ is $C_2$-cofinite and a $\Z_+$-graded weak module $U$ is not 
an Artin module, then there is an irreducible $V$-module $W$ such 
that $\dim \Hom_V(W, U)=\infty$. 
\end{lmm}

\pr
Since $V$ is $C_2$-cofinite, there is a finite set 
$B_{2}(V)\subseteq V$ of homogeneous elements such that 
$C_2(V)+\langle B_{2}(V)\rangle=V$. We first note that $U$ is 
$\Z_+$-graded 
$\oplus_{m=0}^{\infty}U(m)$ by Remark 1(i). Then $U(m)$ is an 
$A_{1,m}(V)$-module, (see Def. 6.2 for the definition of $m$-th Zhu 
algebra $A_{1,m}(V)$). Since $V$ has only finitely many inequivalent 
classes of irreducible modules, the decomposition of $U$ into the 
direct sum of indecomposable $V$-modules is determined by the 
decomposition of $N$-th component $U(N)$ as an $A_{1,N}(V)$-module 
for $N$ sufficiently large.

Consider a submodule 
$U(w)=\langle v_1(i_1)\cdots v_k(i_k)w\mid v_i\in V, i_j\in \Z
\rangle$ of $U$ 
generated from one element $w\in U(N)$. Then from Remark 1(v) we can 
choose a spanning set 
$$\{v_1(n_1)....v_k(n_k)w \mid v_i\in B_{2}(V), 
n_1<\cdots <n_k <N+q, \} $$
where $q$ is the largest weight of elements in $B_{2}(V)$. Therefore 
the character of $U(w)=\sum_{m=0}^{\infty}\dim U(w)(m)e^{2\pi im}$ 
has a upper bound, which implies that the length of composition 
series is bounded, say $K$. Since this holds for any $w\in U(n)$ and 
$\sum_{w\in U(N)} U(w)=U$, The Loewy length of $U$ is less than or 
equal to $K$. Since $V$ is $C_2$-cofinite, $N$-th Zhu algebra 
$A_{1,N}(V)$ is a finite dimensional $\C$-algebra and so 
for any two irreducible $A_{1,N}(V)$-modules $T^1$ and $T^2$, 
the possibility of extension $0\to T^1\to T \to T^2\to 0$ is finite. 
Therefore, if $U$ is not an Artin module, then $\dim{\rm soc}(U)
=\infty$, but $V$ has only finitely many inequivalent simple 
modules. Therefore, we have the desired conclusion. 
\prend

\begin{thm}
If $V$ is $C_2$-cofinite, then $D((W_1\otimes W_2)^{\ast})$ is an 
Artin module for any two Artin modules $W_1$ and $W_2$.  
In particular, a tensor product $W_1\boxtimes_{\la} W_2$ always 
exists as an Artin module and 
$$W_1\boxtimes_{\la} W_2=(D((W_1\otimes_{\la} W_2)^{\ast}))'.$$
The above construction is depend on $\la$, but 
it is uniquely determined as an isomorphism classes of $V$-modules. 
\end{thm}

\pr  
As we showed, the $C_2$-cofiniteness on $V$ implies that a weak 
module $D((W_1\otimes W_2)^{\ast})$ is a $\Z_+$-graded module and is 
also a direct sum of generalized eigenspaces of $L(0)$. 
If $D((W_1\otimes W_2)^{\ast})$ is not an Artin module, then 
there is an irreducible $V$-module $W_3$ such that 
$\dim\Hom_V(W_3, D((W_1\otimes W_2)^{\ast})=\infty$, that is, 
a fusion rule $N_{W_1,W_2}^{(W_3)'}$ is infinite, which contradicts 
Theorem 3.3. Therefore $D((W_1\otimes_{\la}W_2)^{\ast})$ and its 
restricted dual are Artin modules. Denote the restricted dual 
$(D((W_1\otimes_{\la}W_2)^{\ast}))'$ by $S$ for a while. 

There is a natural map $F:W_1\otimes W_2 \to 
(D((W_1\otimes_{\la}W_2)^{\ast})^{\ast}=\ol{S}$. We write 
$w_1\boxtimes_{\la}w_2$ to denote the image of $w_1\otimes w_2$ 
by this map. Then $F$ is clearly a $P(\la)$-intertwining map. We 
next show that $S$ is spanned by all homogenous components of 
elements of the form $\{L(0)^m w_1\boxtimes_{\la}w_2 \mid 
w_1\in W_1, w_2\in W_2, m=0,1,... \}$.
Suppose false, then there is a nonzero element 
$\mu \in D((W_1\otimes_{\la} W_2)^{\ast})$ such that 
$\langle \mu, w_1\boxtimes_{\la}w_2\rangle =0 
\mbox{  for all }w_1\in W_1, w_2\in W_2$. On the other hand, there 
are $w_1\in W_1$ and $w_2\in W_2$ such that 
$\langle \mu,w_1\otimes w_2\rangle \not=0$.  
However, in $(D((W_1\otimes_{\la}W_2)^{\ast}))'$, $w_1\otimes w_2$ 
is equal to $w_1\boxtimes_{\la}w_2$, which contradicts the above. 
Therefore, $S$ is the $\la$-tensor product $W_1\boxtimes_{\la}W_2$ 
of $W_1$ and $W_2$. 
\prend

We have already shown the following in the proof of 
the previous theorem. 

\begin{thm}
$W_1\boxtimes_{\la} W_2$ is spanned by the homogenous elements of 
$w_1\boxtimes_{\la} w_2$ as a $\C[L(0)]$-module. \end{thm}

\section{Associativity of products of intertwining operators}
In this section we will show the associativity and commutativity of 
tensor products. Namely, we will prove the following theorem,  
(cf. Theorem 14.8 in \cite{H1}, Theorem 3.3 in \cite{H2}).

\begin{thm}
Let $V$ be a vertex operator algebra satisfying $C_2$-cofiniteness. 
Then logarithmic intertwining operators for $V$ have the 
following associativity properties: \\
(1) For any Artin modules $W_1, W_2, W_3, W_0'$ and $W_5$ and any 
logarithmic intertwining operators $\CY_1$ and $\CY_2$ of type 
$\binom{W_0'}{W_1W_5}$ and $\binom{W_5}{W_2W_3}$ respectively, 
there exist an Aritin module $W_6$ and logarithmic intertwining 
operators $\CY_3$ and $\CY_4$ of type $\binom{W_6}{W_1W_2}$ and 
$\binom{W_0'}{W_6W_3}$ respectively, such that for any 
$z_1,z_2\in \C$ satisfying $|z_1|>|z_2|>|z_1-z_2|>0$,  
$$\langle w_0, \CY_1(w_1,z_1)\CY_2(w_2,z_2)w_3\rangle =\langle 
w_0, \CY_4(\CY_3(w_1,z_1-z_2)w_2,z_2)w_3\rangle  \eqno{(5.1)}$$
holds for any $w_0\in W_0, w_1\in W_1$, $w_2\in W_2$ and 
$w_3\in W_3$. \\
(2) Conversely, if $\CY_3$, $\CY_4$ of the above types are given, 
then there are an Artin module $W_5$ and logarithmic intertwining 
operators $\CY_1$, $\CY_2$ of the above types satisfying (5.1).\\
\end{thm}

Since the proof for (2) is almost the same as (1), we will prove 
only (1). Set 
$$\begin{array}{l}
\mu_{w_0,w_3}^{z_1,z_2}(w_1,w_2)=\langle w_0, \CY_1(w_1,z_1)
\CY_2(w_2,z_2)w_3\rangle \cr
=\langle w_0, \displaystyle{\sum_{i=0}^{K'}\sum_{j=0}^{H'}\sum_{n,m\in \Q}}(w_1)_{(n,i)}
z_1^{-n-1}\log^i(z_1)(w_2)_{(m,j)}z_2^{-m-1}\log^j(z_2)w_3\rangle
\end{array} \eqno{(5.2)}$$
for $w_i\in W_i$. Viewing $\mu_{w_0,w_3}^{z_1,z_2}(w_1,w_2)$ as a 
function on $z_1$ and $z_2$, it is a solution of differential 
equations of regular singular points. From the theory of 
differential equations (for example, see Appendix in \cite{K}), 
we know that (5.2) is absolutely convergent when $|z_1|>|z_2|>0$ and 
there exist $r_{i,k}$, $s_{i,k}\in \R$ such that (5.2) 
can be analytically extended to the multivalued analytic function 
$$\sum_{i=1}^p\sum_{k=0}^K\sum_{h=0}^H\sum_{m=0}^{\infty}
a_{i,h,k,m}(w_0,w_1,w_2,w_3)z_2^{r_{i}-m}(z_1-z_2)^{s_{i}+m}
\left(\log(z_1-z_2)\right)^k\left(\log(z_2)\right)^h\eqno{(5.3)}$$
in the region $|z_2|>|z_1-z_2|>0$. Since (5.3) for general elements 
$w_i\in W_i$, $i=0,1,2,3$, are 
$\C[z_1^{\pm},z_2^{\pm},(z_1-z_2)^{-1}]$-linear combinations of 
(5.3) for those $w_i\in W_i$ $(i=0,1,2,3)$ satisfying 
$\sum_{i=0}^3\deg(w_i)\leq \sum_{i=0}^3 p_i^1$, we see that $K$ and 
$H$ can be taken to be independent of $w_i\in W_i$ $i=0,1,2,3$. 
We may assume that all $w_0,w_1,w_2,w_3$ are homogeneous. 

Simplify the notation, set $s=\wt(w_0)-\wt(w_1)-\wt(w_2)-\wt(w_3)$. 
First we claim that $r_i, s_i$ are all rational numbers and 
$r_i+s_i=s$. For, we have $-\!n\!-\!1\!-\!m\!-\!1=s$ in (5.2) 
since $\langle w_0, (w_1)_{(n,i)}((w_2)_{(m,j)})w_3)\rangle =0$ 
if $\wt(w_0)\not=\wt((w_1)_{(n,i)}((w_2)_{(m,j)})w_3)$. 
On the other hand, we can rewrite (5.2) into (5.3) by using 
$\log(z_1)=\log(z_2+z_1-z_2)=\log(z_2)+\log(1+\frac{z_1-z_2}{z_2}))
=\log(z_2)+\sum_{i=0}^{\infty}\frac{1}{i}(\frac{-z_1+z_2}{z_2})^i$ 
and $(z_1/z_2)^m=(1+\frac{z_1-z_2}{z_2})^m=
\sum_{i=0}^{\infty}\binom{m}{i}(z_1-z_2)^iz_2^{-i}$. 
Therefore, there is a rational number $t=t(w_0,w_1,w_2,3_3)$ such that 
$\mu_{w_0,w_3}^{z_1,z_2}(w_1,w_2)$ has the form 
$$\sum_{m<t}\sum_{k=0}^K\sum_{h=0}^H
a_{h,k,m}(w_0,w_1,w_2,w_3)z_2^{m}(z_1-z_2)^{s-m}
\left(\log(z_1-z_2)\right)^k\left(\log(z_2)\right)^h.\eqno{(5.5)}$$ 

\begin{lmm}
There is an integer $N$ which does not depend on $w_0,w_1,w_2,w_3$ such that 
$$a_{h,k,m}(w_0,w_1,w_2,w_3)=0 \mbox{  for  }m<N+2\deg(w_0).$$ 
\end{lmm}

\pr
We can choose $N$ so that 
$$a_{h,k,m}(u_0,u_1,u_2,u_3)=0 \mbox{   for  }m<N+2\deg(u_0)$$ 
for all $(u_0,u_1,u_2,u_3)$ satisfying 
$\sum_{i=0}^3 \deg(u_i)\leq p=\sum_{i=0}^3p_i$. 
Since $V$ is $C_1$-cofinite, we have the desired result from 
(3.1) by induction. 
For example, suppose that $w_0\otimes v(-1)w_1\otimes w_2\otimes w_3$ 
is a minimal counterexample. Then from (3.1), we have 
$$ \begin{array}{l}
S^{1-}(w_0\otimes v(-1)w_1\otimes w_2\otimes w_3;z_1,z_2)\cr
\mbox{}\quad =\sum_{k=0}^{\infty}(z_2)^k
S^{1-}(v^{\ast}(-1-k)w_0\otimes w_1\otimes w_2\otimes w_3;z_1,z_2)\cr
\mbox{}\qquad +\sum_{k=0}^{\infty}\binom{-1}{k}(z_1+z_2)^{-1-k}
S^{1-}(w_0\otimes w_1\otimes v(k)w_2\otimes w_3;z_1,z_2) \cr 
\mbox{}\qquad +\sum_{k=0}^{\infty}(z_1)^{-1-k}
S^{1-}(w_0\otimes w_1\otimes w_2\otimes v(k)w_3;z_1,z_2) 
\end{array}$$
Replacing $(z_1)^k$ and $(z_1+z_2)^k$ by 
$(z_2+z_1-z_2)^k$ and $(2z_2+z_1-z_2)^k$, respectively, 
we have that the maximal power of $z_2$ in 
$S^{1-}(w_0\otimes v(-1)w_1\otimes w_2\otimes w_3;z_1,z_2)$ is less 
than $N+2\wt(w_0)$, since $\wt(v^{\ast}(-1-k)w_0)=\wt(w_0)-\wt(v)-k$. 
\prend

Let $\la_1$ and $\la_2$ be complex numbers satisfying 
$|\la_1|>|\la_2|>|\la_1-\la_2|$. 
Keeping $\CY_3(w_1,\la_1-\la_2)w_2
=\sum_{m,h} (w_1)_{(m,h)}w_2(\la_1-\la_2)^{-m-1}\log(\la_1-\la_2)^h$ in mind, 
set 
$$W=\oplus_{h=0}^{\infty}(\oplus_{m\in \C}W_1\otimes_{m,h}W_2))$$ 
and define the action of $v(n)$ on $W$ by 
$$v(n)(w_1\otimes_{m,h}w_2)=w_1\otimes_{m,h}v(n)w_2
+\sum_{i=0}^{\infty}
\binom{n}{i}(\la_1-\la_2)^{n-i}(v(i)w_1)\otimes_{m+n-i,h}w_2. $$
Using (5.5), we define linear maps 
$$ \CY_3\in \Hom(W_1\otimes W_2, \hat{W}) 
\qquad \mbox{  and} \qquad
\CY_4\in \Hom(W, \Hom(W_3,W_0')) $$
by 
$$ \begin{array}{l}
\CY_3(w_1,\la_1-\la_2)w_2=\sum_{m,h}w_1\otimes_{m,h}w_2(\la_1-\la_2)^{-m-1}\log(\la_1-\la_2)^h 
\in \overline{W} \qquad \mbox{  and} \cr 
\langle w_0,\CY_4(w_1\otimes_{m,h}w_2,\la_2)w_3\rangle 
=\sum_{k=0}^Ka_{h,k,m+1+s}(w_0,w_1,w_2,w_3)\la_2^{m+1+s}\log^k(z\la_2).
\end{array}$$
Set $W_6=W/\{w\in W \mid \CY_4(w,\la_2)=0\}$.  
We may view $\CY_3(w_1,\la_1-\la_2)w_2\in \overline{W_6}$ and   
$\CY_4(w,\la_2)$ is well-defined for $w\in W_6$. 
We will show that $W_6$ is a $V$-module and $\CY_3$ and $\CY_4$ are 
intertwining maps. The essential part of the proof is to prove 
lower truncation property of vertex operators on $W_6$, 
because the others easily come from (5.1).
For example, the definition of the action $v(n)$ 
on $\oplus_{m\in \C}(W_1\otimes_{m,h}W_2)$ was given to 
lead $\CY_3$ satisfying the commutativity: 
$$
[v(m),\CY_3(w_1,\la_1-\la_2)]w_2=\CY_3(
\sum_{i=0}^{\infty}
\binom{m}{i}(\la_1-\la_2)^{m-i}v(i)w_1,\la_1-\la_2)w_2.$$ 
Also the following calculation implies the associativity of $\CY_4$. 
$$\begin{array}{rl}
0=&\langle w_0,\CY_1(w_1,\la_1)\CY_2(v(m)w_2,\la_2)w_3\rangle 
-\langle w_0,\CY_4(\CY_3(w_1,\la_1-\la_2)v(m)w_2,\la_2)w_3\rangle\cr
=&\langle w_0,\CY_1(w_1,\la_1)\sum_{i=0}^{\infty}\binom{m}{i}
(-\la_2)^iv(m-i)\CY_2(w_2,\la_2)w_3 \cr
&-\langle w_0,\sum_{i=0}^{\infty}\binom{m}{i}(-1)^{i+m}
\la_2^{m-i}\CY_1(w_1,\la_1)\CY_2(w_2,\la_2)v(i)w_3\rangle  \cr
&-\langle w_0,\CY_4(\CY_3(w_1,\la_1-\la_2)v(m)w_2,\la_2)w_3\rangle
\cr
=&\sum_{i=0}^{\infty}\binom{m}{i}(-\la_2)^i
\langle w_0,v(m-i)\CY_1(w_1,\la_1)\CY_2(w_2,\la_2)w_3\rangle  \cr
&-\sum_{i=0}^{\infty}\binom{m}{i}(-\la_2)^i\sum_{j=0}^{\infty}
\binom{m-i}{j}\la_1^{m-i-j}\langle w_0,\CY_1(v(j)w,\la_1)
\CY_2(w_2,\la_2)w_3\rangle  \cr
&-\sum_{i=0}^{\infty}\binom{m}{i}(-1)^{i+m}\la_2^{m-i}
\langle w_0,\CY_1(w_1,\la_1)\CY_2(w_2,\la_2)v(i)w_3\rangle  \cr
&-\langle w_0,\CY_4(\CY_3(w_1,\la_1-\la_2)v(m)w_2,\la_2)w_3\rangle
\cr
=&\sum_{i=0}^{\infty}\binom{m}{i}(-\la_2)^i\langle w_0,
v(m-i)\CY_4(\CY_3(w_1,\la_1-\la_2)w_2,\la_2)w_3\rangle \cr 
&-\sum_{j=0}^{\infty}(\la_1-\la_2)^{m-j}
\langle w_0,\CY_4(\CY_3(v(j)w,\la_1-\la_2)w_2,\la_2)w_3\rangle \cr
&-\sum_{i=0}^{\infty}\binom{m}{i}(-1)^{i+m}\la_2^{m-i}
\langle w_0,\CY_4(\CY_3(w_1,\la_1-\la_2)w_2,\la_2)v(i)w_3\rangle \cr
&-\langle w_0,\CY_4(v(m)\CY_3(w_1,\la_1-\la_2)w_2,\la_2)w_3\rangle\cr
&+\sum_{j=0}^{\infty}(\la_1-\la_2)^{m-j}
\langle w_0,\CY_4(\CY_3(v(j)w_1,\la_1-\la_2)w_2,\la_2)w_3\rangle \cr
=&\sum_{i=0}^{\infty}\binom{m}{i}(-\la_2)^(m-i)
\langle w_0,v(m-i)\CY_4(\CY_3(w_1,\la_1-\la_2)w_2,\la_2)w_3\rangle\cr
&-\sum_{i=0}^{\infty}\binom{m}{i}(-1)^{i+m}\la_2^{m-i}
\langle w_0,\CY_4(\CY_3(w_1,\la_1-\la_2)w_2,\la_2)v(i)w_3\rangle \cr
&-\langle w_0,\CY_4(v(m)\CY_3(w_1,\la_1-\la_2)w_2,\la_2)w_3\rangle. 
\end{array}$$
Starting from 
$\langle w_0,\CY_1(v(m)w_1,\la_1)\CY(w_2,\la_2)w_3\rangle $, we can 
prove the associativity of $\CY_3$ similarly. The other conditions 
are also easy to check. 

For $m\in \C$, $w_0\in W_0$ and $w_3\in W_3$, define 
$\be_{-m}(w_0,w_3)\in (W_6)^{\ast}$ by 
$$ \begin{array}{l}
 \be_{-m}(w_0,w_3)(w_1\otimes_{n,h}w_2) \cr
\mbox{}\quad =\sum_{k=0}^Ka_{h,k,n,-m}(w_0,w_1,w_2,w_3)
(\la_1-\la_2)^n(\log(\la_1-\la_2))^h \log(\la_1)^k\in \C
\end{array} $$
for $n\in \Q$, $h\in \Z$, $w_1\in W_1$ and $w_2\in W_2$. Then by the definition, we have 
$$\mu_{w_0,w_3}^{\la_1,\la_2}=\sum_{-m\leq N+2\deg(w_0)}
\be_{-m}(w_0,w_3)\la_2^{-m}\cdot \CY_3 \in (W_1\otimes W_2)^{\ast}. $$

We note that since $Y^6(v,\la_1-\la_2)$ on $W_6$ 
is well-defined and satisfies Jacobi-identity, 
$(Y^6)^{\dagger}(v,\la_1-\la_2)$ on $(W_6)^{\ast}$ 
is also well-defined and satisfies Jacobi-identity. 
In particular, $D(W_6^{\ast})$ is a $\Z_+$-graded 
module since $V$ is $C_2$-cofinite.

\begin{prn}  
$\be_{n}(w_0,w_3)\in D(W_6^{\ast})$. 
\end{prn}
  
\pr 
It suffices to show $v(r)\be_n(w_0,w_3)=0$ for $v\in V$ and 
$r>\!\!>0$. 
By the definition, we have 
$$\begin{array}{l}
(v^{\ast}(-1)(\be_{n}(w_0,w_3)))(w_1,w_2)\cr
\mbox{}\quad =(\be_{n}(w_0,w_3))
((\sum_{i=0}^{\infty}\binom{-1}{i}(\la_1-\la_2)^{-1-i}v(i)w_1,w_2)
+(w_1,a(-1)w_2)). \end{array}$$
Consider 
$$S^1(w_0\otimes \sum_{i=0}^{\infty}\binom{-1}{i}(\la_1-\la_2)^{-1-i}
v(i)w_1\otimes w_2\otimes w_3;z_1,z_2) +
S^1(w_0\otimes w_1\otimes v(-1)w_2\otimes w_3)$$ 
and denote it by $S_a(z_1,z_2)$, which is equal to 
$$S^1(\sum_{i=0}^{\infty}z_2^iS^1(v^{\ast}(-1-i)w_0\otimes w_1
\otimes w_2\otimes w_3;z_1,z_2)+
S(w_0\otimes w_1\otimes w_2\otimes \sum_{i=0}z_2^{-1-i}v(i)w_3)$$
by (3.1). If we take $a=L(-1)^rv$ for $r$ sufficiently large, then 
$\frac{1}{r!}a(-1)=v(-r-1)$ and $a^{\ast}(-1-i)w_0=0$ for all 
$i\geq 0$ in the first term and $a(i)=0$ for $0\leq i<r$ in the 
second term. Therefore, $S(z_1,z_2)$ has an expansion 
$$S_{L(-1)^rv}(z_1,z_2)=\sum_{h=0}^H\sum_{k=0}^K
\sum_{-m<N+2\deg(w_0)-1-r}\alpha_mz_2^{-m}(z_1-z_2)^{s+m}\log(z_2)^k
\log(z_1-z_2)^h $$
in $|z_1|>|z_2|>|z_1-z_2|$ for $r$ sufficiently large. 
In particular, it has no power of $z_2^{n}$. That is, 
$$v(r)\be_n(w_0,w_3)=0 \mbox{   for  }r>\!\!>0. $$  
\prend

Let's go back to the proof of Theorem 5.1. 
Using a $V$-homomoprhism $\phi$: $W_1\otimes_{\la_1-\la_2}W_2 \to W_6$ given by 
$$\phi(w_1\otimes_{\la_1-\la_2}w_2)=
\sum w_1\otimes_{m,h}w_2(\la_1-\la_2)^{-m-1}\log^h(\la_1-\la_2),$$
we have $\phi^{\ast}:W_6^{\ast} \to (W_1\otimes_{\la_1-\la_2})^{\ast}$. 
In particular, $Y'(v,z)$ satisfies the Jacobi identity on ${\rm Im}\phi^{\ast}$.  
Since ${\rm Im}\phi^{\ast}$ contains $\be_m(w_0,w_3)$, 
$$\be_{m}(w_0,w_3)\in D((W_1\otimes_{\la_1-\la_2} W_2)^{\ast})$$
and 
$$  \mu_{w_0,w_3}^{\la_1,\la_2}\in 
\ol{D((W_1\otimes_{\la_1-\la_2}W_2)^{\ast})}=(W_1\boxtimes_{\la_1-\la_2}W_2)^{\ast}. $$
Hence by setting $\mu^{\la_1,\la_2}(w_0\otimes w_3)=\mu_{w_0,w_3}^{\la_1,\la_2}$, 
we obtain 
$$\begin{array}{rl}
\mu^{\la_1,\la_2}\in& \Hom(W_0\otimes W_3, \ol{D((W_1\otimes_{\la_1-\la_2}W_2)^{\ast})} \cr
&\cong \Hom(W_1\boxtimes_{\la_1-\la_2}W_2, (W_0\otimes W_3)^{\ast}) \cr
&\cong \Hom((W_1\boxtimes_{\la_1-\la_2}W_2)\otimes W_3, (W_0)^{\ast}) \cr
&=\Hom((W_1\boxtimes_{\la_1-\la_2}W_2)\otimes W_3, \ol{(W_0)'}). 
\end{array}$$
A direct calculations shows that this is a $P(\la)$-intertwining map 
of type $\binom{W_0'}{W_1\boxtimes_{\la_1-\la_2}W_2, W_3}$ and so  
$W_6$ is a homomorphic image of $W_1\boxtimes_{\la_1-\la_2}W_2$. 

This completes the proof of Theorem 5.1. \\

From the direct calculation, we have:
$$\begin{array}{l}
v(n)(w_1\otimes_{\la_1}(w_2 \otimes_{\la_2} w_3))\cr
\mbox{}\quad =(\sum_{i=0}^{\infty} \binom{n}{i}\la_1^{n-i}v(i)w_1)
\otimes (w_2\otimes_{\la_2}w_3) 
+(w_1\otimes_{\la_1}v(n)(w_2 \otimes_{\la_2} w_3)) \cr
\mbox{}\quad =(\sum_{i=0}^{\infty} \binom{n}{i}\la_1^{n-i}v(i)w_1)
\otimes (w_2\otimes_{\la_2}w_3) \cr
\mbox{}\qquad +(w_1\otimes_{\la_1}(\sum_{i=0}^{\infty}\binom{n}{i}
\la_2^{n-i}v(i)w_2 \otimes_{\la_2} w_3))
+(w_1\otimes_{\la_1}(w_2 \otimes_{\la_2} v(n)w_3)) \cr
\mbox{and}\cr 
v(n)((w_1\otimes_{\la_1-\la_2}w_2) \otimes_{\la_2} w_3) \cr
\mbox{}\quad =(\sum_{i=0}^{\infty} \binom{n}{i}(\la_2)^{n-i}v(i)
(w_1\otimes_{\la_1-\la_2}w_2)\otimes_{\la_2}w_3)
+((w_1\otimes_{\la_1-\la_2}w_2) \otimes_{\la_2} v(n)w_3) \cr
\mbox{}\quad =(\sum_{i=0}^{\infty} \binom{n}{i}(\la_2)^{n-i}
(\sum_{j=0}^i\binom{i}{j}(\la_1-\la_2)^{i-j}v(j)w_1
\otimes_{\la_1-\la_2}w_2)\otimes_{\la_2}w_3) \cr
\mbox{}\qquad +(\sum_{i=0}^{\infty} \binom{n}{i}(\la_2)^{n-i}
(w_1\otimes_{\la_1-\la_2}v(i)w_2) \otimes_{\la_2}w_3)
+((w_1\otimes_{\la_1-\la_2}w_2) \otimes_{\la_2} v(n)w_3) \cr
\mbox{}\quad =(\sum_{i=0}^{\infty}\binom{n}{i}(\la_1)^{n-i}v(i)w_1
\otimes_{\la_1-\la_2}w_2)\otimes_{\la_2}w_3) \cr
\mbox{}\qquad +((w_1\otimes_{\la_1-\la_2}\sum_{i=0}^{\infty}
\binom{n}{i}(\la_2)^{n-i}v(i)w_2) \otimes_{\la_2}w_3)
+((w_1\otimes_{\la_1-\la_2}w_2) \otimes_{\la_2} v(n)w_3).
\end{array}$$ 
Since they have similar forms, we can identify 
$W_1\otimes_{\la_1}(W_2\otimes_{\la_2}W_3)$ and 
$(W_1\otimes_{\la_1-\la_2}W_2)\otimes_{\la_2}W_3$ in 
$W_1\otimes W_2\otimes W_3$ by 
$\omega_1\otimes_{\la_1}(w_2\otimes_{\la_2}w_3) \leftrightarrow 
(w_1\otimes_{\la_1-\la_2}w_2)\otimes_{\la_2}w_3$.  
Similarly, we can identity 
$(W_1\otimes_{\la_1}W_2)\otimes_{\la_2}W_3$ and 
$(W_2\otimes_{-\la_1}W_1)\otimes_{\la_2+\la_1}W_3$.

By the similar argument as in the proof of 
Theorem 5.1, we have: 

\begin{thm}  Assume that $V$ is $C_2$-cofinite and $W_1$, $W_2$ and 
$W_3$ are Artin modules. Then there are two unique isomorphisms 
$$\phi^a: W_1\boxtimes_{\la_1}(W_2\boxtimes_{\la_2}W_3) 
\longrightarrow (W_1\boxtimes_{\la_1-\la_2}W_2)\boxtimes_{\la_2}W_3$$
such that $\phi^a(w_1\boxtimes_{\la_1}(w_2\boxtimes_{\la_2}w_3))
=(w_1\boxtimes_{\la_1-\la_2}w_2)\boxtimes_{\la_2}w_3$ 
for $|\la_1|>|\la_2|>|\la_1-\la_2|>0$ and 
$$\phi^c: W_1\boxtimes_{\la_1}(W_2\boxtimes_{\la_2}W_3) 
\longrightarrow 
(W_2\boxtimes_{-\la_1}W_1)\boxtimes_{\la_2+\la_1}W_3$$ 
such that $\phi^c(w_1\boxtimes_{\la_1}(w_2\boxtimes_{\la_2}w_3))
=(w_2\boxtimes_{-\la_1}w_1)\boxtimes_{\la_2+\la_1}w_3$ 
for $|\la_1|>|\la_2|>|\la_1+\la_2|>0$. 
\end{thm}

\section{Application of tensor product theory}

\subsection{Projective modules}
Before we show an application, we have to study a few things about 
Artin modules. 

\begin{dfn} 
An Artin module $W$ is called {\bf projective} if for any 
epimorphism $\phi:U \to W$, there is a homomorphism $\psi:W \to U$ 
such that $\phi\psi=1_W$. \end{dfn}

The following is easy. 

\begin{lmm}
If $W$ is a direct factor of a projective module, then $W$ is 
projective. \end{lmm}

By Remark 1(v), if a weak module $W$ is generated from one element 
$w$, then $W$ is spanned by 
$$\{u_1(i_1)u_2(i_2)\cdot u_m(i_m)w \mid u_j\in B_2(V), i_1<i_2<\cdots <i_m \}. \eqno{(6.1)}$$ As an application of this spanning set, 
we prove the following proposition.

\begin{prn}
For any irreducible module $U$, there is a projective 
module $W$ with a homomorphic image $U$. 
\end{prn}

\pr
Let $\mu$ be a conformal weight of $U$ and consider all weak modules 
generated from one element with weight $\mu$. By the above spanning 
set, such modules are all Artin module and their characters have a 
upper bound. So there is a largest Artin module $W$ among such 
modules with a homomorphic image $U$. We will show that $W$ is 
projective. Let $u\in W$ be an element with $\wt(u)=\mu$ such that 
$W$ is generated from $u$. If there is an epimorphism $\phi:Q\to W$, 
then there is a homogeneous element $u'\in Q$ such that 
$\phi(u')=u$. Then a submodule $W'$ generated from $u'$ is 
isomorphic to $W$ and so $Q=\Ker\phi\oplus W'$. 
\prend

Since an Artin module is generated from a finite set, similar 
arguments show: 

\begin{prn}
For any Artin module $U$, there is a projective module $P$ with a 
homomorphic image $U$.
\end{prn}

\subsection{Orbifold models}
In this subsection, we will consider a simple vertex operator 
algebra. 

Recently, there are several reports about the construction of VOAs 
from small VOAs. In the most cases, using a finite Abelian group 
$G$, a small VOA $W$ and its modules $W^g$ labeled by $g\in G$, 
they construct a $G$-graded VOA $V=\oplus_{g\in G} W^g$. 
Here "$G$-graded" means that 
$\{v(n)w \mid v\in W^g, w\in W^h, n\in \Z\} \subseteq W^{g+h}$ 
for any $g,h\in G$. In this case, $G^{\ast}=\Hom(G,\C)$ becomes an 
automorphism group of $V$ and $W$ is the fixed point subVOA 
$V^{G^{\ast}}$. In the most cases, they showed the rationality and 
$C_2$-cofiniteness using the rationality and $C_2$-cofiniteness of 
$W$. 

Generalizing these setting, we consider the following situation. \\
Let $V$ be a simple VOA and $G$ a finite automorphism group of $V$ 
(including non-Abelian groups) and assume that $V^G$ is 
$C_2$-cofinite. The purpose of this section is to investigate the 
representations of orbifold models $V^{\langle g\rangle }$ for 
$g\in G$. In order to find a $V^{\langle g\rangle }$-module, the 
first step is to find it in 
$V$-modules. The second class is the following modules. 

\begin{dfn}[$g$-twisted module]
Let $g$ be an automorphism of $V$ of order $p$. 
An Artin $g$-twised $V$-module 
$M=\oplus_{n\in \N/p}M(n)$ ($\dim M(n)<\infty$) 
is an $\frac{\Z}{p}$-graded vector space 
equipped with a linear map 
$$ \begin{array}{rcl}
  V &\longrightarrow &(\End(M))\{z\} \cr
  v &\rightarrow  &Y_M(v,z)=\sum_{n\in \frac{1}{p}\Z}v(n)z^{-n-1} 
\end{array}$$
which satisfies the following for all $0\leq r\leq p-1$, 
$u\in V^r=\{u\in V \mid g(u)=e^{2\pi \sqrt{-1}r/p}u\}$, $v\in V$, 
$w\in M$, 
$$\begin{array}{l}
  Y_M(u,z)=\sum_{n\in r/p+\Z}u(n)z^{-n-1}  \cr
  u(m)M(n)\subseteq M(n+\wt(u)-m-1) \cr
  u(m)w=0 \mbox{  for  }m>\!\!>0 \cr
  Y_M(\1,z)=1;
\end{array}$$ 
and the twisted Jacobi form: 
$$\begin{array}{l}
z_0^{-1}\delta\left(\frac{z_1-z_2}{z_0}\right)Y^M(u,z_1)Y^M(v,z_2)
-z_0^{-1}\delta\left(\frac{z_2-z_1}{-z_0}\right)
Y^M(v,z_2)Y^M(u,z_1)\cr
=z_2^{-1}\left(\frac{z_1-z_0}{z_2}\right)^{-r/p}
\delta\left(\frac{z_1-z_0}{z_2}\right)Y^M(Y(u,x_0)v,x_2). 
\end{array}$$
\end{dfn}

Different from ordinary $g$-twisted modules, 
we don't assume the semisimplicity of $L(0)$, (but we won't 
use this difference here). 

We say that $V$ is {\it $g$-rational} if every Artin $g$-twisted 
$V$-module is completely reducible. 

In order to investigate the structure of $g$-twisted modules, the 
most important tool is a $g$-twisted $n$-th Zhu algebra $A_{g,n}(V)$ 
of $V$ associated with $g\in {\rm Aut}(V)$ of order $p$ and 
nonnegative $n\in \frac{1}{p}{\mathbb Z}$. 
This algebra was first introduced in \cite{Z} for the case when 
$g$ is the identity element and $n=0$ and then 
Dong, Li and Mason extended it to general case in \cite{DLM}. 

\begin{dfn} 
Fix nonnegative $n=l+i/p\in {\mathbb Z}/p$ with $l$ a nonnegative 
integer and $0\leq i\leq p-1$. For $0\leq r\leq p-1$ we define 
$\delta_i(r)=1$ if $i\geq r$ and $\delta_i(r)=0$ if $i<r$.
We also set $\delta_i(p)=1$. Set 
$$O_{g,n}(V)=\left< 
\begin{array}{l}
 (L(-1)+L(0))v, \cr 
\mbox{Res}_zY(u,z)v\frac{(1+z)^{\wt(u)-1+\delta_i(r)+l+r/p}}
{z^{2l+\delta_i(r)+\delta_i(p-r)}} \mid 
u,v\in V, g(u)=e^{2\pi \sqrt{-1}r/p}u. 
\end{array}\right> $$
and $A_{g,n}(V)=V/O_{g,n}(V)$. The product $*_{g,n}$ on $V$ for 
homogeneous $u\in V$ and $v\in V$ is defined by 
$$ u*_{g,n}v=\frac{1}{T}\sum_{m=0}^{l}(-1)^m\binom{m+l}{l}
\mbox{Res}_zY(\frac{1}{p}\sum_{i=0}^{p-1} g^i(u),z)v
\frac{(1+z)^{\wt(u)+l}}{z^{l+m+1}}$$
Extend this linearly to obtain a bilinear product on $V$.
\end{dfn}

One of important properties of $g$-twisted $n$-th Zhu algebra is 
that if $U=\oplus_{m\in \Z/p}U(m)$ is a $g$-twisted module, then 
$U(m)$ is an $A_{g,m}(V)$-module and conversely 
if $T$ is an $A_{g,n}(V)$-module, then there is a $g$-twisted 
module $U=\oplus_{n\in \N/p}U(n)$ such that $U(n)\cong T$ as 
$A_{g,n}(V)$-modules, see \cite{DLM}. Therefore, $A_{g,n}(V)$ are 
semisimple algebras for all $n$ if and only if $V$ is $g$-rational.

As in the above definition, a product in $A_{g,n}$ depends on $g$. 
However, the existence of uniformed product was shown in \cite{MT}. 

\begin{prn}[\cite{MT}]
Let $V$ be a simple VOA and $G$ a finite automorphism group of $V$. 
Set $$A^G_{n}(V)=V/\cap_{g\in G}O_{g,n}(V).$$
Then there is a uniform product in $A_n^G(V)$ such that 
$$\begin{array}{ccc}
A^G_n(V)&\longrightarrow  &\oplus_{g\in G}A_{g,n}(V) \cr
  v     &\to               &(v,...,v) 
\end{array}$$
is a ring isomorphism.  
\end{prn}

By the definition of $O_{g,n}(V)$, we have the following by the same 
argument as in \cite{GN} (cf. Th. 2.4 in \cite{Miy}). 

\begin{lmm}  If $V^G$ is $C_2$-cofinite, then 
$A_{g,n}(V)$ is finite dimensional for any $g\in G$ and 
$n\in \N/|G|$. 
\end{lmm}

We now apply the tensor product theory to obtain the 
following theorem.

\begin{thm}
Let $V$ be a simple vertex operator algebra and 
$g$ an automorphism of $V$ of finite order. If a fixed point 
subVOA $V^{\langle g\rangle }$ is $C_2$-cofinite, then every 
irreducible $V^{\langle g\rangle }$-module is contained in a 
$h$-twisted $V$-module for some $h\in \langle g\rangle $. 
\end{thm}

\pr 
Let $p$ be an order of $g$ and set 
$$V^{n}=\{v\in V \mid g(v)=e^{2\pi \sqrt{-1}n/p}v \} 
\mbox{  for   }n\in \Z/p\Z. $$
        Let $U$ be an irreducible $V^{0}$-module and $W$ a 
projective indecomposable module with a homomorphic image $U'$. 
We will show that there are an integer $k$ and an $g^k$-twisted 
module $T$ containing $W$ as a sub $V^{\langle g\rangle }$-module. 
Then $T'$ is the desired one. 

By using the products in $V$, there are natural epimorphisms 
$$ \phi:\underbrace{V^1\boxtimes \cdots \boxtimes V^1}_p \to V^0  
\qquad \mbox{and} \qquad  
\phi_W:\underbrace{V^1\boxtimes \cdots \boxtimes V^1}_p\boxtimes  
W\to W.$$
Since $W$ is projective, $\underbrace{V^1\boxtimes \cdots 
\boxtimes V^1}_p\boxtimes  W$ has a direct factor 
$\tilde{W}\cong W$ and so for each $i=0,1,...,p$, 
$\underbrace{V^1\boxtimes \cdots \boxtimes V^1}_i\boxtimes  W$ has 
a direct factor $Q_i$ such that $Q_{i+1}$ is a direct factor of 
$V^1\boxtimes Q_i$ and $Q_p=\tilde{W}$. If we define $S$ by 
$$\begin{array}{l}
S(u^1,u^2,...,u^{p};z_1,...,z_{p},z)\cr
\mbox{}\qquad =\langle w_0,Y^W(Y(\cdots (Y(Y(u^1,z_1-z_2)u^2,z_2-z_3)
\cdots, z_{p-1}-z_p)u_p,z_p)w\rangle  
\end{array} \eqno{(6.2)}$$
for $u^1,\cdots,u^p\in V^1$, $w\in W$ and $w_0\in W'$, 
then it is absolutely convergent when 
$0<|z_1-z_2|<|z_2-z_3|<\cdots <|z_{p-1}-z_p|<|z_p|$ and we 
extend it to the whole space except for singular points 
by analytic continuations. By the choices of $Q_i$, 
there are (logarithmic) intertwining operators $\CY_i$ of type 
$\binom{Q_{i+1}}{V^1,Q^i}$ such that 
$$S(w_0,u^1,u^2,...,u^{p},w;z_1,...,z_{p})=
\langle w_0, \CY_1(u^1,z_1)\cdots \CY_p(u^p,z_p)w\rangle  $$
in the region $|z_1|>\cdots >|z_p|>0$. 
Since (6.2) has no $(\log(z_i))$-terms, $\CY_i(u^i,z_i)$ has no 
$(\log(z_i))$-term for every $i$. Using the products of intertwining 
maps given by tensor products, we can express $S$ by 
$$ \langle w_0, u^1\boxtimes_{z_1}(u^2\boxtimes_{z_2} \cdots 
\boxtimes_{z_{p-1}} (u^p\boxtimes_{z_p} w)\cdots)\rangle . $$
Then by the associativity of tensor products, we have 
$$\begin{array}{l}
\langle w_0, ((\cdots (u^1\boxtimes_{z_1-z_2} u^2)\boxtimes \cdots  
\boxtimes_{z_{p-1}-z_p} u^p)\boxtimes_{z_p} w\rangle \cr
\mbox{}\qquad \qquad \qquad =\langle w_0, u^1
\boxtimes_{z_1}(u^2\boxtimes_{z_2} \cdots 
\boxtimes_{z_{p-1}}(u^p\boxtimes_{z_p} w)\cdots)\rangle  \cr 
\langle w_0, w_1\boxtimes_{z_{i-1}-z_{i+2}} 
((u^i\boxtimes_{z_i-z_{i+1}} u^{i+1})
\boxtimes_{z_{i+1}-z_{i+2}} w_2)\rangle  \cr
\mbox{}\qquad \qquad \qquad =
\langle w_0, w_1\boxtimes_{z_{i-1}-z_{i+2}} 
(u^i\boxtimes_{z_{i}-z_{i+2}} 
(u^{i+1}\boxtimes_{z_{i+1}-z_{i+2}} w_2))\rangle 
\end{array} \eqno{(6.3)}$$
for $w_1=((\cdots (u^1\boxtimes_{z_1-z_2} u^2)
\boxtimes_{z_2-z_3} ... \boxtimes_{z_{i-2}-z_{i-1}} u^{i-1})$ and 
$w_2=(u^{i+2}\boxtimes_{z_{i+2}}( .... \boxtimes_{z_{p-1}} 
(u^p\boxtimes w)\cdots)$. We note that tensor product 
$u^i\boxtimes_{z_i-z_{i+1}}u^{i+1}$ is induced by the products 
in $V$.

Before we continue the proof, we need the following lemma. 

\begin{lmm}  Let $\CY^1$, $\CY^2$ and $\CY^3$ be 
intertwining operators of 
types $\binom{Q_{i+2}}{V^1,Q_{i+1}}$, $\binom{Q_{i+1}}{V^1,Q_i}$ 
and $\binom{Q_{i+2}}{V^2,Q_i}$. Assume that  
$$\langle w_0,\CY^1(u,z_1)\CY^2(v,z_2)w\rangle =\langle w_1,
\CY^3(Y(u,z_1-z_2)v,z_2)w\rangle  $$
hold for all $u,v\in V^1$ and $w_1\in Q_{i+2}'$, $w\in Q_i$ 
in the region $|z_1|>|z_2|>|z_1-z_2|>0$. Then 
$$  \langle w_1,\CY^3(Y(u,z_1-z_2)v,z_2)w\rangle =\langle w_1,
\CY^1(v,z_2)\CY^2(u,z_1)w\rangle  $$
always hold in $|z_2|>|z_1-z_2|>|z_1|>0$.  
In particular, as functions extended by analytic continuations, 
we have 
$$  \langle w_1,\CY^1(u,z_1)\CY^2(v,z_2)w\rangle =\langle w_1,
\CY^1(v,z_2)\CY^2(u,z_1)w\rangle . \eqno{(6.4)}$$
In particular, if $Q_i,Q_{i+1}$ are indecomposable, then the powers 
of $z$ in $\CY^1(v,z)$ and those of $z$ in $\CY^2(v,z)$ are 
congruent modulo integers. 
\end{lmm}

\pr Using skew-symmetry, we have the following equations. 
$$\begin{array}{l}
\langle w_0,\CY^1(u,z_1)\CY^2(v,z_2)w\rangle \cr
=\langle w_0,\CY^3(Y(u,z_1-z_2)v,z_2)w\rangle  
=\langle w_0,\CY^3(e^{(z_1-z_2)L(-1)}Y(v,-z_1+z_2)u,z_2)w\rangle\cr
=\langle w_0,\CY^3(Y(v,-z_1+z_2)u, z_2+z_1-z_2)w\rangle  
=\langle w_0,\CY^1(v,z_2)\CY^2(u,z_2+z_1-z_2)w\rangle . 
\end{array}$$
Therefore, we have the first desired result (6.4). 
Since 
$S(w_0,w;z_1,z_2)$ \\
$(=\langle w_0,\CY_1(u,z_1)\CY_2(v,z_2)w\rangle )$ 
satisfies differential equations of regular singular points with 
singular points $z_1=0$, $z_2=0$ and $z_1-z_2=0$ at most and 
$S(w_0,w_1;z_1,z_2)$ has no $\log(z_1)$, $\log(z_2)$, 
$\log(z_1-z_2)$-terms, $S(w_0,w_1;z_1,z_2)$ is a finite sum of 
forms $\frac{f_j(z_1,z_2)}{z_1^{a_j}z_2^{b_j}(z_1-z_2)^{c_j}}$ in 
the complex plane, where $f_j(z_1,z_2)$ is holomorphic on $\C^2$. 
Then each has a formal power series of types 
$$\begin{array}{lcl}
z_1^{-a_j-c_j-m}z_2^{-b_j+m} &in &\langle w_0,
\CY_1(u,z_1)\CY_2(v,z_2)w\rangle   \cr
z_1^{-a_j+m}z_2^{-b_j-c_j-m} &in &\langle w_0,
\CY_1(v,z_2)\CY_2(u,z_1)w\rangle   \cr
z_2^{-b_j-a_j-m}(z_1-z_2)^{-c_j+m} &in &\langle w_0,
\CY_3(Y(u,z_1-z_2)v,z_2)w\rangle  
\end{array}$$
with $m\in \N$. Clearly, $c_j\in \Z$. We note that since 
$Q_i, Q_{i+1}, Q_{i+2}$ are indecomposable and $V^1$ is irreducible, 
the powers of $z_i$ in $\CY_i$ are uniquely determined up to 
integers, that is, $a_i\equiv a_j \pmod{\Z}$ and 
$b_i\equiv b_j \pmod{\Z}$ for $i,j$. Therefore, $a_i\equiv b_j 
\mod{\Z}$ for $i,j$ and we have the desired result.  
\prend

It follows from Lemma 6.10 that the powers of $z_i$ in 
$\CY_i(u,z_i)$ are congruent each other modulo integers . Since the 
total sums are integers, there is an integer $m$ such that the 
powers of $z_i$ in $\CY_i(u,z_i)$ are $m/p+\Z$. 

Using the above information, let's construct a twisted module 
$T=W\oplus Q_1\oplus \cdots \oplus Q_{p-1}$.  

Firstly, for $u\in V^0\oplus V^1$, we define a vertex operator 
$Y^T(u,z)$ by a matrix 
$$ Y^T(u,z)=\begin{pmatrix}T_{11}(u,z)&\cdots &T_{1p}(u,z) \cr 
    \vdots &\cdots & \vdots \cr
 T_{p1}(u,z)&\cdots & T_{pp}(u,z) \end{pmatrix}, $$
where $T_{ij}$ are defined by 
$$ \left\{\begin{array}{l}
 T_{ii}(v,z)=Y^{Q_i}(v,z) \mbox{  and  }T_{ij}(v,z)=0 
\mbox{  if  }i\not=j \qquad \mbox{for }u\in V^0 \cr
T_{i,i+1}(u,z)=\CY_i(u,z)\mbox{  and  }T_{ij}(u,z)=0 
\mbox{  if  }j\not\equiv i+1 \pmod{p} \qquad \mbox{for }u\in V^1. 
\end{array}\right.$$
Then for any $t\in T$, the powers of $z$ in $Y^T(u,z)t$ have a lower 
bound. It is also easily from Lemma 6.10 that we have 
$$  \langle t_0, Y^T(u,z_1)Y^T(u',z_2)t\rangle =\langle t_0, 
Y^T(u',z_2)Y^T(u,z_1)t\rangle  $$
for all $t\in T$ and $t_0\in T'$. Therefore, all vertex operators of 
type $Y^T(u,z)$ $(u\in V^1)$ and $Y^T(v,z)$ $v\in V^0$ are locally 
commutative each others. Clearly they satisfy $L(-1)$-derivative 
property. By the associativity of products of intertwining operators 
and $\underbrace{V^1\boxtimes \cdots \boxtimes V^1}_p\boxtimes W 
\to V^0\boxtimes W \to W$, we can define vertex operators of all 
elements in $V$. Therefore if $(m,p)=1$, then choosing $k\in \Z$ 
such that $hk\equiv 1 \pmod{p}$, $T$ is a $g^k$-twisted module. 
If $(m,p)=n\not=1$, then $T$ is a $g^{nk}$-twisted module for some $k$. 
\prend

\begin{thm}
Let $V$ be a simple vertex operator algebra and 
$G$ a finite automorphism group of $V$. If a fixed point subVOA 
$V^G$ is rational and $C_2$-finite, then $V$ is $g$-rational and 
$V^{\langle g\rangle }$ is rational for any $g\in G$.
\end{thm}

\pr
Since $V^G$ is $C_2$-cofinite and $V^{\langle g\rangle }$ is a 
direct sum of finite number of irreducible $V^G$-module by 
\cite{DM}, $V^{\langle g\rangle }$ is also $C_2$-cofinite for any 
$g\in G$. As we mentioned in \S 2, since 
$A_{1,n}(V^{\langle g\rangle })$ is a finite dimensional algebra 
and $A_{g,n}(V)$ is a homomorphic image of 
$A_{1,l}(V^{\langle g\rangle })$, $A_{g,n}(V)$ is also a finite 
dimensional algebra, where $l$ is the maximal integer not exceed $n$.

For $g,h\in G$ and a $g$-twisted module $W$, 
$hW$ becomes a $hgh^{-1}$-twisted module. 
Furthermore, if  $W$ has a $H$-action for 
some $H\leq C_G(g)$, then 
$\oplus_{h\in G/H} hW$ is a direct sum of twisted modules with 
a $G$-action on it. We call such a module twisted $GV$-module. 

As the first step, we will show that if $W$ is a nonzero twisted 
$GV$-module, then $W^G\not=\{0\}$. As they showed in \cite{DM}, $V$ 
contains all irreducible $G$-modules as a $G$-module. 
Let $0\not=U\subseteq W$ be an irreducible $G$-submodule of $W$ and 
$T\subseteq V$ be a $G$-submodule of $V$ isomorphic to $\Hom(U,\C)$. 
Then as they showed in \cite{HMT}, using the results in \cite{DM}, 
we can prove that $W$ contains a $G$-submodule isomorphic to 
$U\otimes T$. In particular, $W$ contains a trivial $G$-submodule.

Let $A^G_n(V)=V/\cap_{g\in G}O_{g,n}(V)$ be a $G$-twisted $n$-th Zhu 
algebra introduced in \cite{MT}. Since $(A^G_n(V))^G$ is a 
homomorphic image of $n$-th Zhu algebra $A_{1,n}(V^G)$ for a fixed 
point subVOA $V^G$, $(A^G_n(V))^G$ is semisimple by the assumption. 
In particular, it does not contain any nonzero nilpotent ideal. 
Therefore, $(J(A_n^G(V)))^G=\{0\}$ for every $n\in \frac{1}{|G|}\Z$, 
where $J(R)$ denotes the Jacobson radical (the largest nilpotent 
ideal) of a ring $R$.

If $V$ is not $g$-rational for some $g$ then $A_{g,n}(V)$ is not 
semisimple for $n$ sufficiently large and so $J(A^G_n(V))\not=0$ 
since $A_{g,n}(V)$ is a homomorphic image of $A^G_n(V)$.   

Since $V^G$ is $C_2$-cofinite, the number of irreducible 
$V^G$-modules is finite and $t={\rm max}\{m \mid W(m)=0, 
W(m+1)\not=0 \mbox{  for some $V^G$-module }W\}$ is well defined. 
If we take $n\in \Z$ greater than $t$, then 
all modules induced from irreducible $A^G_n(V)$-modules becomes 
irreducible. Let $W$ be a $V$-module induced from 
$A^G_n(V)$-module $A^G_n(V)$. Then $W$ is $G$-invariant and is a sum 
of twisted modules. By the construction, $W(m)$ is isomorphic to 
$A^G_m(V)$ as an $A^G_m(V)$-module for $m>t$. Let $J(W)$ be a 
submodule induced from $A_n^G(V)$-submodule $J(A^G_n(V))$  
($\subseteq A^G_n(V)$).
$J(W)$ is also $G$-invariant and is a sum of twisted modules.  
Similarly, $J(W)(m)$ is isomorphic to $J(A^G_m(V))$ for $m\geq n$. 
However, since $J(W)$ is $G$-invariant, $J(W)^G\not=0$. Since 
$J(W)^G$ is a $V^G$-module, $(J(W)(m))^G=J(A^G_m(V))^G\not=0$ for 
$m$ sufficiently large, which contradicts $J(A^G_m(V))\not=0$. 

That is, $V$ is $h$-rational for any $h\in G$. By the previous 
theorem, every irreducible $V^{\langle g\rangle }$-module is 
contained in a $h$-twisted $V$-module, all 
$V^{\langle g\rangle }$-modules are completely reducible. 
Namely, $V^{\langle g\rangle }$ is rational. 
This completes the proof of Theorem 6.11. 
\prend

As an application of these theorems, we can easily determine the 
representations of an orbifold model $(V_L)^+=(V_L)^{\theta}$ of 
lattice VOA $V_L$ given by an automorphism $\theta$ of order two 
induced from $-1$ on $L$, where $L$ is a general $n$-dimensional 
positive definite even lattice since we have already known that 
$(V_{\Z a})^{+}$ of rank one satisfies $C_2$-finite condition by 
\cite{Y}. Namely, find a sublattice $\Z a_1\perp \cdots \perp 
\Z a_n \subseteq L$. Then there is an automorphism group $G$ 
such that $(V_L)^G \cong (V_{\Z a_1})^+\otimes \cdots 
\otimes (V_{\Z a_n})^+$. Since $V_L^G$ is $C_2$-cofinite, the above 
theorem gives an alternative proof of the following result, 
which was recently shown in \cite{AD}:  
every irreducible $V_L^+$-module is contained in a $V_L$-module or 
$\theta$-twisted $V$-module. 

\begin{cry}
Every simple module of $(V_L)^{+}$ is a submodule of 
$V_L$-module or a $\theta$-twisted $V_L$-module.  
If $V_{\Z a}^+$ is rational for rank one lattice $\Z a$, 
then $V_L^+$ is also rational. 
\end{cry}

\end{document}